\documentclass[oneside,english]{amsart}
\usepackage[T1]{fontenc}
\usepackage[latin9]{inputenc}
\usepackage{amsthm}
\usepackage{amssymb}
\usepackage{setspace}
\usepackage{esint}
\makeatletter
\numberwithin{equation}{section}
\numberwithin{figure}{section}
\theoremstyle{plain}
\newtheorem{thm}{\protect\theoremname}
  \theoremstyle{plain}
  \newtheorem{cor}[thm]{\protect\corollaryname}
  \theoremstyle{remark}
  \newtheorem*{rem*}{\protect\remarkname}
  \theoremstyle{plain}
  \newtheorem{lem}[thm]{\protect\lemmaname}

\makeatother

\usepackage{babel}
  \providecommand{\corollaryname}{Corollary}
  \providecommand{\lemmaname}{Lemma}
  \providecommand{\remarkname}{Remark}
\providecommand{\theoremname}{Theorem}

\begin{document}

\title{On the Hausdorff and packing measures of slices of dynamically defined
sets}

\author{Ariel Rapaport}

\date{February 18, 2015}

\subjclass[2000]{Primary: 28A80 Secondary: 28A78}

\keywords{Self-similar set, Self-affine set, Hausdorff measures, packing measures}

\thanks{Supported by ERC grant 306494}
\begin{abstract}
Let $1\le m<n$ be integers, and let $K\subset\mathbb{R}^{n}$ be
a self-similar set satisfying the strong separation condition, and
with $\dim K=s>m$. We study the a.s. values of the $s-m$-dimensional
Hausdorff and packing measures of $K\cap V$, where $V$ is a typical
$n-m$-dimensional affine subspace. For $0<\rho<\frac{1}{2}$ let
$C_{\rho}\subset[0,1]$ be the attractor of the IFS $\{f_{\rho,1},f_{\rho,2}\}$,
where $f_{\rho,1}(t)=\rho\cdot t$ and $f_{\rho,2}(t)=\rho\cdot t+1-\rho$
for each $t\in\mathbb{R}$. We show that for certain numbers $0<a,b<\frac{1}{2}$,
for instance $a=\frac{1}{4}$ and $b=\frac{1}{3}$, if $K=C_{a}\times C_{b}$
then typically we have $\mathcal{H}^{s-m}(K\cap V)=0$.
\end{abstract}
\maketitle

\section{Introduction}

Let $1\le m<n$ be integers, and given $0\le t\le n$ let $\mathcal{H}^{t}$
and $\mathcal{P}^{t}$ be the $t$-dimensional Hausdorff and Packing
measures respectively. Let $s\in(m,n)$ be a real number, and let
$K\subset\mathbb{R}^{n}$ be compact with $0<\mathcal{H}^{s}(K)<\infty$.
Denote by $\mu$ the restriction of $\mathcal{H}^{s}$ to $K$, by
$G$ the set of all $n-m$-dimensional linear subspaces of $\mathbb{R}^{n}$,
and by $\xi_{G}$ the natural measure on $G$. It is well known that
$\dim_{H}(K\cap(x+V))=s-m$ and $\mathcal{H}^{s-m}(K\cap(x+V))<\infty$,
for $\mu\times\xi_{G}$-a.e. $(x,V)\in K\times G$ (see Theorem 10.11
in \cite{key-1}). It is also known that if $s=\dim_{P}K$ then $\dim_{P}(K\cap(x+V))\le\max\{0,s-m\}$
for every $V\in G$ and $\mathcal{H}^{m}$-a.e. $x\in V^{\perp}$
(see Lemma 5 in \cite{key-12}), where $\dim_{P}$ stands for the
packing dimension. In this paper $K$ will denote certain self-similar
or self-affine sets, in which cases it will be shown that more can
be said about the $\mu\times\xi_{G}$-typical values of $\mathcal{H}^{s-m}(K\cap(x+V))$
and $\mathcal{P}^{s-m}(K\cap(x+V))$.

Assume first that $K$ is a self-similar set which satisfies the strong
separation condition (SSC). If $m=1$ and $K$ is rotation-free, then
from a result by Kempton (Theorem 6.1 in \cite{key-6}) it follows
that $\mathcal{H}^{s-m}(K\cap(x+V))>0$ for $\mu\times\xi_{G}$-a.e.
$(x,V)$, if and only if $\frac{dP_{V^{\perp}}\mu}{d\mathcal{H}^{m}}\in L^{\infty}(dP_{V^{\perp}}\mu)$
for $\xi_{G}$-a.e. $V$, where $P_{V^{\perp}}$is the orthogonal
projection onto $V^{\perp}$. In Theorem \ref{T4} below the case
of a general $1\le m<n$ and a general self-similar set $K$, satisfying
the SSC, will be considered. A necessary and sufficient condition
for $\mathcal{H}^{s-m}(K\cap(x+V))>0$ to holds for $\mu\times\xi_{G}$-a.e.
$(x,V)$ will be given. In Corollary \ref{C2} this condition is verified
when $m=1$, $s>2$ and the rotation group of $K$ is finite. Also
given in Theorem \ref{T4}, is a necessary and sufficient condition
for $\mathcal{H}^{s-m}(K\cap(x+V))=0$ to hold for $\mu\times\xi_{G}$-a.e.
$(x,V)$.

Continuing to assume that $K$ is a self-similar set with the SSC,
it will be shown in Theorem \ref{T2} that $\mathcal{P}^{s-m}(K\cap(x+V))>0$
for $\mu\times\xi_{G}$-a.e. $(x,V)$. Also given in Theorem \ref{T2},
is a sufficient condition for $\mathcal{P}^{s-m}(K\cap(x+V))=\infty$
to hold for $\mu\times\xi_{G}$-a.e. $(x,V)$. By using this condition,
it is shown in Corollary \ref{C1} that this is in fact the case when
$m=1$ and $s>2$. This extends a result of Orponen (Theorem 1.1 in
\cite{key-5}), which deals with the case in which $n=2$, $s>m=1$
and $K$ is rotation-free.

Lastly we consider the case in which $n=2$, $m=1$ and $K$ is a
certain self-affine set. For $0<\rho<\frac{1}{2}$ let $C_{\rho}\subset[0,1]$
be the attractor of the IFS $\{f_{\rho,1},f_{\rho,2}\}$, where $f_{\rho,1}(t)=\rho\cdot t$
and $f_{\rho,2}(t)=\rho\cdot t+1-\rho$ for each $t\in\mathbb{R}$.
It will be assumed that $K=C_{a}\times C_{b}$, where $0<a,b<\frac{1}{2}$
are such that $a^{-1}$ and $b^{-1}$ are Pisot numbers, $\frac{\log b}{\log a}$
is irrational, and $\dim_{H}(C_{a})+\dim_{H}(C_{b})>1$. Under these
conditions it is shown in \cite{key-7} that there exists a dense
$G_{\delta}$ set, of $1$-dimensional linear subspaces $V\subset\mathbb{R}^{2}$,
such that $P_{V}\mu$ and $\mathcal{H}^{1}$ are singular. By using
this fact, it will be proven in Theorem \ref{T6} below that $\mathcal{H}^{s-m}(K\cap(x+V))=0$
for $\mu\times\xi_{G}$-a.e. $(x,V)$. This result demonstrates some
kind of smallness of the slices $K\cap(x+V)$, hence it may be seen
as related to a conjecture made by Furstenberg (Conjecture 5 in \cite{key-13}).
In our setting this conjecture basically says that for $\xi_{G}$-a.e.
$V\in G$ we have $\dim_{H}(K\cap(x+V))\le\max\{\dim_{H}K-1,0\}$
for each $x\in\mathbb{R}^{2}$, which demonstrates the smallness of
the slices in another manner.

The rest of this article is organized as follows: In section \ref{S2}
the results are stated. In section \ref{S3} the results regarding
self-similar sets are proven. In section \ref{S4} we prove the aforementioned
theorem regarding self-affine sets.

\textbf{Acknowledgment. }I would like to thank my adviser Michael
Hochman, for suggesting to me problems studied in this paper, and
for many helpful discussions.

\section{\label{S2}statement of the results}

\subsection{Slices of self-similar sets}

Let $0<m<n$ be integers, let $G$ be the Grassmann manifold consisting
of all $n-m$-dimensional linear subspaces of $\mathbb{R}^{n}$, let
$O(n)$ be the orthogonal group of $\mathbb{R}^{n}$, and let $\xi_{O}$
be the Haar measure corresponding to $O(n)$. Fix $U\in G$ and for
each Borel set $E\subset G$ define 
\begin{equation}
\xi_{G}(E)=\xi_{O}\{g\in O(n)\:|\: gU\in E\},\label{E15}
\end{equation}
then $\xi_{G}$ is the unique rotation invariant Radon probability
measure on $G$. For a linear subspace $V$ of $\mathbb{R}^{n}$ let
$P_{V}$ be the orthogonal projection onto $V$, let $V^{\perp}$
be the orthogonal complement of $V$, and set $V_{x}=x+V$ for each
$x\in\mathbb{R}^{n}$.

Let $\Lambda$ be a finite and nonempty set. Let $\{\varphi_{\lambda}\}_{\lambda\in\Lambda}$
be a self-similar IFS in $\mathbb{R}^{n}$, with attractor $K\subset\mathbb{R}^{n}$
and with $\dim_{H}K=s>m$. For each $\lambda\in\Lambda$ there exist
$0<r_{\lambda}<1$, $h_{\lambda}\in O(n)$ and $a_{\lambda}\in\mathbb{R}^{n}$,
such that $\varphi_{\lambda}(x)=r_{\lambda}\cdot h_{\lambda}(x)+a_{\lambda}$
for each $x\in\mathbb{R}^{n}$. We assume that $\{\varphi_{\lambda}\}_{\lambda\in\Lambda}$
satisfies the strong separation condition. Let $H$ be the smallest
closed sub-group of $O(n)$ which contains $\{h_{\lambda}\}_{\lambda\in\Lambda}$,
and let $\xi_{H}$ be the Haar measure corresponding to $H$. For
each $E\subset\mathbb{R}^{n}$ set $\mu(E)=\frac{\mathcal{H}^{s}(K\cap E)}{\mathcal{H}^{s}(K)}$,
then $\mu$ is a Radon probability measure which is supported on $K$.

For each $0\le s<\infty$, $\nu$ a Radon probability measure on $\mathbb{R}^{n}$,
and $x\in\mathbb{R}^{n}$ set
\begin{equation}
\Theta^{*s}(\nu,x)=\underset{\epsilon\downarrow0}{\limsup}\:\frac{\nu(B(x,\epsilon))}{(2\epsilon)^{s}}\mbox{\:\ and }\:\Theta_{*}^{s}(\nu,x)=\underset{\epsilon\downarrow0}{\liminf}\:\frac{\nu(B(x,\epsilon))}{(2\epsilon)^{s}},\label{E13}
\end{equation}
where $B(x,\epsilon)$ is the closed ball in $\mathbb{R}^{n}$ with
center $x$ and radios $\epsilon$. It holds that $\Theta^{*s}(\nu,\cdot)$
and $\Theta_{*}^{s}(\nu,\cdot)$ are Borel functions (see remark 2.10
in \cite{key-1}). For $V\in G$ define
\[
F_{V}(x,h)=\Theta_{*}^{m}(P_{(hV)^{\perp}}\mu,P_{(hV)^{\perp}}(x))\:\mbox{ for }(x,h)\in K\times H,
\]
then $F_{V}$ is a Borel function from $K\times H$ to $[0,\infty]$.
In what follows the collection $\{F_{V}\}_{V\in G}$ will be of great
importance for us.

Let $\mathcal{V}$ be the set of all $V\in G$ with
\[
\xi_{H}(H\setminus\{h\in H\::\: P_{(hV)^{\perp}}\mu\ll\mathcal{H}^{m}\})=0\:.
\]
In Lemma \ref{L1} below it will be shown that $\xi_{G}(G\setminus\mathcal{V})=0$.
First we state our results regarding the Hausdorff measure of typical
slices of $K$.
\begin{thm}
\label{T4}(i) Given $V\in\mathcal{V}$, if $\left\Vert F_{V}\right\Vert _{L^{\infty}(\mu\times\xi_{H})}<\infty$
then $\mathcal{H}^{s-m}(K\cap(x+hV))>0$ for $\mu\times\xi_{H}$-a.e.
$(x,h)\in K\times H$.

(ii) Given $V\in\mathcal{V}$, if $\left\Vert F_{V}\right\Vert _{L^{\infty}(\mu\times\xi_{H})}=\infty$
then $\mathcal{H}^{s-m}(K\cap(x+hV))=0$ for $\mu\times\xi_{H}$-a.e.
$(x,h)\in K\times H$.

(iii) $\mathcal{H}^{s-m}(K\cap V_{x})>0$ for $\mu\times\xi_{G}$-a.e.
$(x,V)\in K\times G$ if and only if $\left\Vert F_{V}\right\Vert _{L^{\infty}(\mu\times\xi_{H})}<\infty$
for $\xi_{G}$-a.e. $V\in G$.

(iv) $\mathcal{H}^{s-m}(K\cap V_{x})=0$ for $\mu\times\xi_{G}$-a.e.
$(x,V)\in K\times G$ if and only if $\left\Vert F_{V}\right\Vert _{L^{\infty}(\mu\times\xi_{H})}=\infty$
for $\xi_{G}$-a.e. $V\in G$.
\end{thm}

From Theorem \ref{T4} we can derive the following corollaries.
\begin{cor}
\label{C2}Assume $m=1$, $s>2$ and $|H|<\infty$, then $\mathcal{H}^{s-m}(K\cap V_{x})>0$
for $\mu\times\xi_{G}$-a.e. $(x,V)\in K\times G$.
\end{cor}

\begin{cor}
\label{C3}Assume that $H=O(n)$ and
\[
\mu\times\xi_{G}\{(x,V)\in K\times G\::\:\mathcal{H}^{s-m}(K\cap V_{x})>0\}>0,
\]
then there exists $0<M<\infty$ such that for each $V\in G$ we have
$P_{V^{\perp}}\mu\ll\mathcal{H}^{m}$ with $\left\Vert \frac{dP_{V^{\perp}}\mu}{d\mathcal{H}^{m}}\right\Vert _{L^{\infty}(\mathcal{H}^{m})}\le M$.
\end{cor}

\begin{rem*}
It is known that under the assumptions of Corollary \ref{C3} we have
$\dim(P_{V^{\perp}}\mu)=m$ for each $V\in G$ (see Theorem 1.6 in
\cite{key-10}). It is not known however if $P_{V^{\perp}}\mu\ll\mathcal{H}^{m}$
for each $V\in G$, which is in fact a major open problem. Hence Corollary
\ref{C3} implies that determining whether
\[
\mu\times\xi_{G}\{(x,V)\in K\times G\::\:\mathcal{H}^{s-m}(K\cap V_{x})>0\}>0
\]
is probably quite hard.
\end{rem*}

Next we state our results regarding the packing measure of typical
slices.
\begin{thm}
\label{T2}(i) $\mathcal{P}^{s-m}(K\cap V_{x})>0$ for $\mu\times\xi_{G}$-a.e.
$(x,V)\in K\times G$.

(ii) Given $V\in\mathcal{V}$, if $\left\Vert \frac{1}{F_{V}}\right\Vert _{L^{\infty}(\mu\times\xi_{H})}=\infty$
then $\mathcal{P}^{s-m}(K\cap(x+hV))=\infty$ for $\mu\times\xi_{H}$-a.e.
$(x,h)\in K\times H$.

(iii) If $\left\Vert \frac{1}{F_{V}}\right\Vert _{L^{\infty}(\mu\times\xi_{H})}=\infty$
for $\xi_{G}$-a.e. $V\in G$, then $\mathcal{P}^{s-m}(K\cap V_{x})=\infty$
for $\mu\times\xi_{G}$-a.e. $(x,V)\in K\times G$.
\end{thm}

From Theorem \ref{T2} the following corollary can be derived.
\begin{cor}
\label{C1}Assume $m=1$ and $s>2$, then $\mathcal{P}^{s-m}(K\cap V_{x})=\infty$
for $\mu\times\xi_{G}$-a.e. $(x,V)\in K\times G$.
\end{cor}

\begin{rem*}
In the proofs of Corollaries \ref{C2} and \ref{C1}, we use the fact
that if $m=1$ and $s>2$ then $\frac{dP_{V^{\perp}}\mu}{d\mathcal{H}^{m}}$
is a continuous function for $\xi_{G}$-a.e. $V\in G$ (see Lemma
3.2 in \cite{key-9} and the discussion before it). It is not known
whether this is still true if $m>1$ or $m<s\le2$, hence we need
the assumptions $m=1$ and $s>2$. 
\end{rem*}

\subsection{Slices of self-affine sets}

Assume $n=2$ and $m=1$. Given $0<\rho<\frac{1}{2}$ define $f_{\rho,1},f_{\rho,2}:\mathbb{R}\rightarrow\mathbb{R}$
by
\[
f_{\rho,1}(x)=\rho\cdot x\:\mbox{ and }\: f_{\rho,2}(x)=\rho\cdot x+1-\rho\quad\mbox{ for each }x\in\mathbb{R},
\]
let $C_{\rho}\subset[0,1]$ be the attractor of the IFS $\{f_{\rho,1},f_{\rho,2}\}$,
set $d_{\rho}=\dim_{H}C_{\rho}$ (so that $d_{\rho}=\frac{\log2}{\log\rho^{-1}}$),
and for each $E\subset\mathbb{R}$ set $\mu_{\rho}(E)=\frac{\mathcal{H}^{d_{\rho}}(C_{\rho}\cap E)}{\mathcal{H}^{d_{\rho}}(C_{\rho})}$.
\begin{thm}
\label{T6}Let $0<a<b<\frac{1}{2}$ be such that $\frac{1}{a}$ and
$\frac{1}{b}$ are Pisot numbers, $\frac{\log b}{\log a}$ is irrational
and $d_{a}+d_{b}>1$, then $\mathcal{H}^{d_{a}+d_{b}-1}((C_{a}\times C_{b})\cap V_{(x,y)})=0$
for $\mu_{a}\times\mu_{b}\times\xi_{G}$-a.e. $(x,y,V)\in C_{a}\times C_{b}\times G$.
\end{thm}

\begin{rem*}
Recall that every integer greater than $1$ is a Pisot number, hence
Theorem \ref{T6} applies for instance in the case $a=\frac{1}{4}$
and $b=\frac{1}{3}$.
\end{rem*}

\begin{rem*}
Note that $0<\mathcal{H}^{d_{a}+d_{b}}(C_{a}\times C_{b})<\infty$,
see Lemma \ref{L10} below.
\end{rem*}

\section{\label{S3}Proof of the results on self similar sets}

\subsection{\label{S3.1}Preliminaries}

The following notations will be used in the proofs of theorems \ref{T4}
and \ref{T2}. For each $\lambda\in\Lambda$ set $p_{\lambda}=r_{\lambda}^{s}$.
Then $\mu$ is the unique self-similar probability measure corresponding
to the IFS $\{\varphi_{\lambda}\}_{\lambda\in\Lambda}$ and the probability
vector $(p_{\lambda})_{\lambda\in\Lambda}$, i.e. $\mu$ satisfies
the relation $\mu=\sum_{\lambda\in\Lambda}p_{\lambda}\cdot\mu\circ\varphi_{\lambda}^{-1}$.
Given a word $\lambda_{1}\cdot...\cdot\lambda_{l}=w\in\Lambda^{*}$
we write $p_{w}=p_{\lambda_{1}}\cdot...\cdot p_{\lambda_{l}}$, $r_{w}=r_{\lambda_{1}}\cdot...\cdot r_{\lambda_{l}}$
, $h_{w}=h_{\lambda_{1}}\cdot...\cdot h_{\lambda_{l}}$, $\varphi_{w}=\varphi_{\lambda_{1}}\circ...\circ\varphi_{\lambda_{l}}$
and $K{}_{w}=\varphi_{w}(K)$. For each $l\geq1$ and $x\in K$, let
$w_{l}(x)\in\Lambda^{l}$ be the unique word of length $l$ which
satisfies $x\in K_{w_{l}(x)}$. Set also 
\begin{equation}
\rho=\min\{d(\varphi_{\lambda_{1}}(K),\varphi_{\lambda_{2}}(K))\::\:\lambda_{1},\lambda_{2}\in\Lambda\mbox{ and }\lambda_{1}\ne\lambda_{2}\},\label{E14}
\end{equation}
then $\rho>0$ since $\{\varphi_{\lambda}\}_{\lambda\in\Lambda}$
satisfies the strong separation condition. Given $V_{1},V_{2}\in G$
set $d_{G}(V_{1},V_{2})=\left\Vert P_{V_{1}}-P_{V_{2}}\right\Vert $
(where $\left\Vert \cdot\right\Vert $ stands for operator norm),
then $d_{G}$ is a metric on $G$.

The following dynamical system will be used in the proofs of theorems
\ref{T4} and \ref{T2}. Set $X=K\times H$ and for each $(x,h)\in X$
let $T(x,h)=(\varphi_{w_{1}(x)}^{-1}x,h_{w_{1}(x)}^{-1}\cdot h)$.
It is easy to check that the system $(X,\mu\times\xi_{H},T)$ is measure
preserving, and from corollary 4.5 in \cite{key-3} it follows that
it is ergodic. Also, for $k\geq1$ and $(x,h)\in X$ it is easy to
verify that $T^{k}(x,h)=(\varphi_{w_{k}(x)}^{-1}x,h_{w_{k}(x)}^{-1}\cdot h)$.

Let $\mathcal{R}$ be the Borel $\sigma$-algebra of $\mathbb{R}^{n}$.
For each $V\in G$ set $\mathcal{R}_{V}=P_{V^{\perp}}^{-1}(\mathcal{R})$,
and let $\{\mu_{V,x}\}_{x\in\mathbb{R}^{n}}$ be the disintegration
of $\mu$ with respect to $\mathcal{R}_{V}$ (see section 3 of \cite{key-8}).
For $\mu$-a.e. $x\in\mathbb{R}^{n}$ the probability measure $\mu_{V,x}$
is defined and supported on $K\cap V_{x}$. Also, for each $f\in L^{1}(\mu)$
the map that takes $x\in\mathbb{R}^{n}$ to $\int f\: d\mu_{V,x}$
is $\mathcal{R}_{V}$-measurable, the formula
\[
\int f\: d\mu=\int\int f(y)\: d\mu_{V,x}(y)\: d\mu(x)
\]
is satisfied, and for $\mu$-a.e. $x\in V^{\perp}$ we have
\[
\int f\: d\mu_{V,x}=\underset{\epsilon\downarrow0}{\lim}\:\frac{1}{P_{V^{\perp}}\mu(B(x,\epsilon))}\cdot\underset{P_{V^{\perp}}^{-1}(B(x,\epsilon))}{\int}f\: d\mu\:.
\]
For more details on the measures $\{\mu_{V,x}\}_{x\in\mathbb{R}^{n}}$
see section 3 of \cite{key-8} and the references therein.

\subsection{Auxiliary lemmas}

We shall now prove some lemmas that will be needed later on. The following
lemma will be used with $\xi_{H}$ in place of $\eta$, when $\xi_{H}$
is considered as a measure on $O(n)$ (which is supported on $H$).
\begin{lem}
\label{L2}Let $Q$ be a compact metric group, and let $\nu$ be its
normalized Haar measure. Let $\eta$ be a Borel probability measure
on $Q$, then for each Borel set $E\subset Q$
\[
\nu(E)=\int_{Q}\eta(E\cdot q^{-1})\: d\nu(q)\:.
\]

\end{lem}

\emph{Proof:} For each Borel set $E\subset Q$ define $\zeta(E)=\int_{Q}\eta(E\cdot q^{-1})\: d\nu(q)$.
Since $\nu$ is invariant it follows that for each $g\in Q$
\[
\zeta(Eg)=\int_{Q}\eta(E\cdot g\cdot q^{-1})\: d\nu(q)=\int_{Q}\eta(E\cdot g\cdot(q\cdot g)^{-1})\: d\nu(q)=\zeta(E)\:.
\]
This shows that $\zeta$ is a right-invariant Borel Probability measure
on $Q$, hence $\nu=\zeta$ by the uniqueness of the Haar measure,
and the lemma follows. $\square$
\begin{lem}
\label{L1}Let $\mathcal{V}$ be the set of all $V\in G$ with
\[
\xi_{H}(H\setminus\{h\in H\::\: P_{(hV)^{\perp}}\mu\ll\mathcal{H}^{m}\})=0,
\]

then $\xi_{G}(G\setminus\mathcal{V})=0$.
\end{lem}

\emph{Proof of Lemma \ref{L1}:} Set $L=G\setminus\{V\in G\::\: P_{V^{\perp}}\mu\ll\mathcal{H}^{m}\}$.
Since $s>m$ it follows that $I_{m}(\mu)<\infty$ (where $I_{m}(\mu)$
is the $m$-energy of $\mu$), hence from theorem 9.7 and equality
(3.10) in \cite{key-1} we get that $\xi_{G}(L)=0$. Let $U\in G$
be as in (\ref{E15}) and set $L'=\{g\in O(n)\::\: gU\in L\}$, then
$\xi_{O}(L')=\xi_{G}(L)=0$. Let $B\subset O(n)$ be a Borel set with
$L'\subset B$ and $\xi_{O}(B)=0$, then from Lemma \ref{L2} it follows
that
\[
0=\xi_{O}(B)=\int\xi_{H}(B\cdot g^{-1})\: d\xi_{O}(g)\:.
\]
From this we get that for $\xi_{O}$-a.e. $g\in O(n)$ 
\begin{multline*}
0=\xi_{H}(B\cdot g^{-1})\geq\xi_{H}(L'\cdot g^{-1})=\xi_{H}\{h\in H\::\: hg\in L'\}=\\
=\xi_{H}(H\setminus\{h\in H\::\: P_{(hgU)^{\perp}}\mu\ll\mathcal{H}^{m}\})\:,
\end{multline*}
and so
\[
\xi_{H}(H\setminus\{h\in H\::\: P_{(hV)^{\perp}}\mu\ll\mathcal{H}^{m}\})=0\:\mbox{ for \ensuremath{\xi_{G}}-a.e. \ensuremath{V\in G},}
\]
which proves the lemma. $\square$
\begin{lem}
\label{L3}Let $\mathcal{Z}$ be the set of all $(x,V)\in K\times G$
such that \textup{$\mu_{V,x}$ is defined and 
\[
\mu_{V,x}(K_{w})=\underset{\epsilon\downarrow0}{\lim}\frac{\mu(K_{w}\cap P_{V^{\perp}}^{-1}(B(P_{V^{\perp}}x,\epsilon)))}{P_{V^{\perp}}\mu(B(P_{V^{\perp}}x,\epsilon))}\:\mbox{ for each }w\in\Lambda^{*},
\]
}then for each $V\in G$ we have 
\[
\mu\times\xi_{H}\{(x,h)\in X\::\:(x,hV)\notin\mathcal{Z}\}=0\:.
\]

\end{lem}

\emph{Proof:} Fix $V\in G$. It holds that $\mathcal{Z}$ is a Borel
set, see section 3 of \cite{key-2} for a related argument. It follows
that the set
\[
\mathcal{Z}_{V}=\{(x,h)\in X\::\:(x,hV)\in\mathcal{Z}\}
\]
is also a Borel set. From the properties stated in section \ref{S3.1}
we get that 
\[
\mu\{x\in K\::\:(x,h)\notin\mathcal{Z}_{V}\}=0\mbox{ for each }h\in H,
\]
and so $\mu\times\xi_{H}(X\setminus\mathcal{Z}_{V})=0$ by Fubini's
theorem. This proves the lemma. $\square$
\begin{lem}
\label{L4} Given a compact set $\tilde{K}\subset\mathbb{R}^{n}$
and $0<t\leq n$, the map that takes $(x,V)\in\tilde{K}\times G$
to $\mathcal{H}^{t}(\tilde{K}\cap V_{x})$ is Borel measurable.
\end{lem}

\emph{Proof:} For $\delta>0$ let $\mathcal{H}_{\delta}^{t}$ be as
defined in section 4.3 of \cite{key-1}. Let $(x,V)\in\tilde{K}\times G$,
$\epsilon>0$ and $\{(x_{k},V^{k})\}_{k=1}^{\infty}\subset\tilde{K}\times G$,
be such that $(x_{k},V^{k})\overset{k}{\rightarrow}(x,V)$. Let $W_{1},W_{2},...\subset\mathbb{R}^{n}$
be open sets with $\tilde{K}\cap V_{x}\subset\cup_{j=1}^{\infty}W_{j}$,
\[
\sum_{j=1}^{\infty}(diam(W_{j}))^{t}\leq\mathcal{H}_{\delta}^{t}(\tilde{K}\cap V_{x})+\epsilon
\]
and $diam(W_{j})\leq\delta$ for each $j\geq1$. Since $\tilde{K}$
is compact and since $(x_{k},V^{k})\overset{k}{\rightarrow}(x,V)$,
it follows that $\tilde{K}\cap V_{x_{k}}^{k}\subset\cup_{j=1}^{\infty}W_{j}$
for each $k\geq1$ which is large enough, and so for each such $k$
\[
\mathcal{H}_{\delta}^{t}(\tilde{K}\cap V_{x_{k}}^{k})\le\sum_{j=1}^{\infty}(diam(W_{j}))^{t}<\mathcal{H}_{\delta}^{t}(\tilde{K}\cap V_{x})+\epsilon.
\]
It follows that the function that maps $(x,V)$ to $\mathcal{H}_{\delta}^{t}(\tilde{K}\cap V_{x})$
is upper semi-continuous, and so Borel measurable. Now since $\mathcal{H}^{s}=\underset{k\rightarrow\infty}{\lim}\mathcal{H}_{1/k}^{s}$
the lemma follows. $\square$
\begin{lem}
\label{L5}Given $0<t\leq n$ and a Radon probability measure $\nu$
on $K\times G$, the map that takes $(x,V)\in K\times G$ to $\mathcal{P}^{t}(K\cap V_{x})$
is $\nu$-measurable (i.e. this map is universally measurable).
\end{lem}

\emph{Proof:} Let $a\geq0$ and set $E=\{(x,V)\in K\times G\::\:\mathcal{P}^{t}(K\cap V_{x})<a\}$,
then in order to prove the lemma it suffice to show that $E$ is $\nu$-measurable.$\newline$Set
$Y=\{C\subset K\::\: C\mbox{ is compact}\}$, endow $Y$ with the
Hausdorff metric, and let $\mathcal{G}$ be the $\sigma$-algebra
of $Y$ which is generated by its analytic subsets. Set 
\[
\mathcal{E}=\{C\in Y\::\:\mathcal{P}^{t}(C)<a\},
\]
then from Theorem 4.2 in \cite{key-14} it follows that $\mathcal{E}\in\mathcal{G}$,
and so from Theorem 21.10 in \cite{key-4} we get that $\mathcal{E}$
is universally measurable.$\newline$For each $(x,V)\in K\times G$
set $\psi(x,V)=K\cap V_{x}$, it will now be shown that $\psi:K\times G\rightarrow Y$
is a Borel function. For each $y\in K$ the function that maps $(x,V)\in K\times G$
to $d(K\cap V_{x},y)$ is lower semi-continuous, and hence a Borel
function. For each $l\geq1$ let $S_{l}\subset K$ be finite and $l^{-1}$-spanning,
and set $\psi_{l}(x,V)=\{y\in S_{l}\::\: d(K\cap V_{x},y)\leq l^{-1}\}$
for each $(x,V)\in K\times G$. It holds that $\psi_{l}:K\times G\rightarrow Y$
is a Borel function and $\psi_{l}\overset{l\rightarrow\infty}{\longrightarrow}\psi$
pointwise, hence $\psi$ is a Borel function. Note also that $E=\psi^{-1}(\mathcal{E})$.$\newline$Since
$\mathcal{E}$ is universally measurable it is $\nu\circ\psi^{-1}$-measurable,
and so there exist $\mathcal{A}$ and $\mathcal{C}$, Borel subsets
of $Y$, such that $\mathcal{A}\subset\mathcal{E}\subset\mathcal{C}$
and $\nu\circ\psi^{-1}(\mathcal{C}\setminus\mathcal{A})=0$. It holds
that $\psi^{-1}(\mathcal{A})$ and $\psi^{-1}(\mathcal{C})$ are Borel
subsets of $K\times G$, $\psi^{-1}(\mathcal{A})\subset E\subset\psi^{-1}(\mathcal{C})$
and $\nu(\psi^{-1}(\mathcal{C})\setminus\psi^{-1}(\mathcal{A}))=0$.
This shows that $E$ is $\nu$-measurable, and the lemma is proved.
$\square$
\begin{lem}
\label{L6}For $(x,h,V)\in K\times H\times G$ set $\psi(x,h,V)=(x,hV)$
and let $B\in K\times G$ be universally measurable. Assume that for
$\xi_{G}$-a.e. $V\in G$ it holds for $\xi_{H}$-a.e. $h\in H$ that
\[
\mu\{x\in K\::\:\psi(x,h,V)\in B\}=0\:,
\]
then $\mu\times\xi_{G}(B)=0$.
\end{lem}

\emph{Proof:} Since $B$ is universally measurable there exist Borel
sets $A,C\subset K\times G$ with $A\subset B\subset C$ and $\mu\times\xi_{H}\times\xi_{G}(\psi^{-1}(C\setminus A))=0$.
From the assumption on $B$ and from Fubini's theorem it follows that
\begin{multline*}
\mu\times\xi_{H}\times\xi_{G}(\psi^{-1}(C))=\mu\times\xi_{H}\times\xi_{G}(\psi^{-1}(A))=\\
=\int\int\mu\{x\::\:(x,h,V)\in\psi^{-1}(A)\}\: d\xi_{H}(h)\: d\xi_{G}(V)\leq\\
\leq\int\int\mu\{x\::\:(x,h,V)\in\psi^{-1}(B)\}\: d\xi_{H}(h)\: d\xi_{G}(V)=0\:.
\end{multline*}
Now from Fubini's theorem, from the definition of $\xi_{G}$ given
in (\ref{E15}), and from Lemma \ref{L2}, it follows that
\begin{multline*}
0=\mu\times\xi_{H}\times\xi_{G}(\psi^{-1}(C))=\\
=\int\int\xi_{H}\{h\::\:(x,h,V)\in\psi^{-1}(C)\}\: d\xi_{G}(V)\: d\mu(x)=\\
=\int\int\xi_{H}\{h\::\:(x,h,gU)\in\psi^{-1}(C)\}\: d\xi_{O}(g)\: d\mu(x)=\\
=\int\int\xi_{H}\{h\::\:(x,hgU)\in C\}\: d\xi_{O}(g)\: d\mu(x)=\\
=\int\int\xi_{H}(\{h\::\:(x,hU)\in C\}\cdot g^{-1})\: d\xi_{O}(g)\: d\mu(x)=\\
=\int\xi_{O}\{g\::\:(x,gU)\in C\}\: d\mu(x)=\\
=\int\xi_{G}\{V\::\:(x,V)\in C\}\: d\mu(x)=\mu\times\xi_{G}(C)\ge\mu\times\xi_{G}(B),
\end{multline*}
which completes the proof of the lemma. $\square$

\subsection{\label{S3.3}Proofs of Theorems \ref{T4} and \ref{T2}}

Fix $V\in\mathcal{V}$ for the remainder of this section, set $F=F_{V}$,
and for each $h\in H$ set $V^{h}=hV$ and $P_{h}=P_{(V^{h})^{\perp}}$.
Set
\[
Q=\{(x,h)\in X\::\: F(x,h)\ne\Theta^{*m}(P_{h}\mu,P_{h}(x))\mbox{ or }F(x,h)=\infty\mbox{ or }F(x,h)=0\}
\]
where $\Theta^{*m}$ is as defined in (\ref{E13}), then $Q$ is a
Borel set. From theorem 2.12 in \cite{key-1} it follows that 
\[
\mu\{x\in K\::\:(x,h)\in Q\}=0\:\mbox{ for each }h\in H\mbox{ with }P_{h}\mu\ll\mathcal{H}^{m},
\]
hence since $V\in\mathcal{V}$ we have 
\begin{equation}
\mu\times\xi_{H}(Q)=\int_{H}\mu\{x\::\:(x,h)\in Q\}\: d\xi_{H}(h)=0\:.\label{E0}
\end{equation}

Let $D$ be the set of all $(x,h)\in X$ such that $P_{h}\mu\ll\mathcal{H}^{m}$,
$\mu_{V^{h},x}$ is defined,
\[
\mu_{V^{h},x}(K_{w})=\underset{\epsilon\downarrow0}{\lim}\frac{\mu(K_{w}\cap P_{h}^{-1}(B(P_{h}x,\epsilon)))}{P_{h}\mu(B(P_{h}x,\epsilon))}\mbox{ for each }w\in\Lambda^{*},
\]
and
\[
0<F(x,h)=\underset{\epsilon\downarrow0}{\lim}\frac{P_{h}\mu(B(P_{h}(x),\epsilon))}{(2\epsilon)^{m}}<\infty\:.
\]
From the choice of $V$, from Lemma \ref{L3} and from (\ref{E0}),
it follows that $\mu\times\xi_{H}(X\setminus D)=0$. Set $D_{0}=\cap_{j=0}^{\infty}T^{-j}D$,
then $\mu\times\xi_{H}(X\setminus D_{0})=0$ since $T$ is measure
preserving. The following lemma will be used several times below.
\begin{lem}
\label{L7}Given $k\geq1$ and $(x,h)\in D_{0}$, we have
\[
\mu_{V^{h},x}(K_{w_{k}(x)})=\left(F(x,h)\right)^{-1}\cdot r_{w_{k}(x)}^{s-m}\cdot F(T^{k}(x,h))\:.
\]

\end{lem}

\emph{Proof:} Set $u=w_{k}(x)$, then
\begin{multline*}
\mu_{V^{h},x}(K_{u})=\underset{\epsilon\downarrow0}{\lim}\:\frac{\mu(K_{u}\cap P_{h}^{-1}(B(P_{h}x,\epsilon)))}{P_{h}\mu(B(P_{h}x,\epsilon))}=\\
=\underset{\epsilon\downarrow0}{\lim}\:\frac{(2\epsilon)^{m}}{P_{h}\mu(B(P_{h}x,\epsilon))}\cdot\frac{\mu(K_{u}\cap P_{h}^{-1}(B(P_{h}x,\epsilon)))}{(2\epsilon)^{m}}=\\
=\left(F(x,h)\right)^{-1}\cdot\underset{\epsilon\downarrow0}{\lim}\:\frac{\mu(K_{u}\cap P_{h}^{-1}(B(P_{h}x,\epsilon)))}{(2\epsilon)^{m}}\:.
\end{multline*}
For each $\epsilon>0$ set $E_{\epsilon}=P_{h_{u}^{-1}h}^{-1}(B(P_{h_{u}^{-1}h}(\varphi_{u}^{-1}(x)),\epsilon\cdot r_{u}^{-1}))$,
then since 
\begin{multline*}
P_{h}^{-1}(B(P_{h}x,\epsilon))=x+V^{h}+B(0,\epsilon)=\\
=\varphi_{u}(\varphi_{u}^{-1}(x)+V^{h_{u}^{-1}h}+B(0,\epsilon\cdot r_{u}^{-1}))=\varphi_{u}(E_{\epsilon}),
\end{multline*}
it follows that
\begin{multline*}
\mu_{V^{h},x}(K_{u})=\left(F(x,h)\right)^{-1}\cdot\underset{\epsilon\downarrow0}{\lim}\:\frac{\mu(\varphi_{u}(K\cap E_{\epsilon}))}{(2\epsilon)^{m}}=\\
=\left(F(x,h)\right)^{-1}\cdot\underset{\epsilon\downarrow0}{\lim}\:\frac{1}{(2\epsilon)^{m}}\sum_{w\in\Lambda^{k}}p_{w}\cdot\mu(\varphi_{w}^{-1}(\varphi_{u}(K\cap E_{\epsilon})))\:.
\end{multline*}
Given $w\in\Lambda^{k}\setminus\{u\}$ we have $\varphi_{u}(K)\cap\varphi_{w}(K)=\emptyset$,
so $\varphi_{w}^{-1}(\varphi_{u}(K))\cap K=\emptyset$, and so
\begin{multline*}
\mu_{V^{h},x}(K_{u})=\left(F(x,h)\right)^{-1}\cdot\underset{\epsilon\downarrow0}{\lim}\:\frac{p_{u}\cdot\mu(K\cap E_{\epsilon})}{(2\epsilon)^{m}}=\left(F(x,h)\right)^{-1}\cdot r_{u}^{s-m}\cdot\underset{\epsilon\downarrow0}{\lim}\:\frac{\mu(E_{\epsilon})}{(2\epsilon\cdot r_{u}^{-1})^{m}}=\\
=\left(F(x,h)\right)^{-1}\cdot r_{u}^{s-m}\cdot F(\varphi_{u}^{-1}(x),h_{u}^{-1}h)=\left(F(x,h)\right)^{-1}\cdot r_{u}^{s-m}\cdot F(T^{k}(x,h)),
\end{multline*}
which proves the lemma. $\square$$\newline$$\newline$\emph{Proof
of theorem \ref{T4}, part (i)}: Assume that $V$ is such that $\left\Vert F\right\Vert _{L^{\infty}(\mu\times\xi_{H})}<\infty$.
Set $M=\left\Vert F\right\Vert _{L^{\infty}(\mu\times\xi_{H})}$ ,
$E=\{(x,h)\::\: F(x,h)\leq M\}$ and $E_{1}=D_{0}\cap(\cap_{j=0}^{\infty}T^{-j}(E))$,
then $\mu\times\xi_{H}(X\setminus E_{1})=0$. For $\xi_{H}$-a.e.
$h\in H$ we have 
\[
\mu\{x\in K\::\:(x,h)\notin E_{1}\}=0,
\]
fix such $h_{0}\in H$. For each $l\geq1$ set
\[
A_{l}=\{x\in K\::\:(x,h_{0})\in E_{1}\mbox{ and }F(x,h_{0})\geq l^{-1}\},
\]
and fix $l_{0}\geq1$. Set $\kappa=\min\{r_{\lambda}\::\:\lambda\in\Lambda\}$,
it will now be shown that
\begin{equation}
\Theta^{*s-m}(\mu_{V^{h_{0}},x},x)\leq(2\rho\kappa)^{m-s}l_{0}M\:\mbox{ for each }x\in A_{l_{0}},\label{E2}
\end{equation}
where $\rho$ is as defined in (\ref{E14}). Let $x\in A_{l_{0}}$
and let $\kappa\rho>\delta>0$. Let $k\geq1$ be such that $r_{w_{k}(x)}\geq\frac{\delta}{\rho}>r_{w_{k+1}(x)}$,
and set $u=w_{k}(x)$. From Lemma \ref{L7} and from $T^{k}(x,h_{0})\in E$
we get that
\[
\mu_{V^{h_{0}},x}(K_{u})=\left(F(x,h_{0})\right)^{-1}\cdot r_{u}^{s-m}\cdot F(T^{k}(x,h_{0}))\leq l_{0}\cdot r_{u}^{s-m}\cdot M\:,
\]
and so
\begin{multline*}
\frac{\mu_{V^{h_{0}},x}(B(x,\delta))}{(2\delta)^{s-m}}\leq\frac{\mu_{V^{h_{0}},x}(B(x,\rho\cdot r_{w_{k}(x)}))}{(2\rho\cdot r_{w_{k+1}(x)})^{s-m}}\leq\\
\leq\frac{\mu_{V^{h_{0}},x}(K_{u})}{(2\rho\kappa\cdot r_{u})^{s-m}}\leq\frac{l_{0}r_{u}^{s-m}M}{(2\rho\kappa\cdot r_{u})^{s-m}}=(2\rho\kappa)^{m-s}l_{0}M,
\end{multline*}
which proves (\ref{E2}).$\newline$It holds that 
\[
\{x\in K\::\:(x,h_{0})\in E_{1}\}=\cup_{l=1}^{\infty}A_{l},
\]
hence
\[
0=\mu(K\setminus\cup_{l=1}^{\infty}A_{l})=\int\mu_{V^{h_{0}},x}(K\setminus\cup_{l=1}^{\infty}A_{l})\: d\mu(x),
\]
and so for $\mu$-a.e. $x\in K$ there exist $l_{x}\geq1$ with $\mu_{V^{h_{0}},x}(A_{l_{x}}\cap V_{x}^{h_{0}})=\mu_{V^{h_{0}},x}(A_{l_{x}})>0$.
Fix such $x_{0}\in K$ and let $y\in A_{l_{x_{0}}}\cap V_{x_{0}}^{h_{0}}$,
then from (\ref{E2}) we get that 
\[
\Theta^{*s-m}(\mu_{V^{h_{0}},x_{0}},y)=\Theta^{*s-m}(\mu_{V^{h_{0}},y},y)\leq(2\rho\kappa)^{m-s}l_{x_{0}}M,
\]
and so from Theorem 6.9 in \cite{key-1} it follows that 
\begin{multline*}
\mathcal{H}^{s-m}(K\cap V_{x_{0}}^{h_{0}})\geq\mathcal{H}^{s-m}(A_{l_{x_{0}}}\cap V_{x_{0}}^{h_{0}})\geq\\
\geq2^{-(s-m)}(2\rho\kappa)^{s-m}l_{x_{0}}^{-1}M^{-1}\cdot\mu_{V^{h_{0}},x_{0}}(A_{l_{x_{0}}}\cap V_{x_{0}}^{h_{0}})>0.
\end{multline*}
This proves that if $\left\Vert F_{V}\right\Vert _{L^{\infty}(\mu\times\xi_{H})}<\infty$,
then for $\xi_{H}$-a.e. $h\in H$ we have 
\[
\mathcal{H}^{s-m}(K\cap(x+hV))>0\:\mbox{ for }\mu\mbox{-a.e. }x\in K,
\]
and so \emph{(i)} follows from Lemma \ref{L4} and Fubini's theorem.\emph{$\newline$$\newline$Proof
of part (ii):} Assume that $V$ is such that $\left\Vert F\right\Vert _{L^{\infty}(\mu\times\xi_{H})}=\infty$,
then 
\[
\mu\times\xi_{H}\{(x,h)\::\: F(x,h)>M\}>0\:\mbox{ for each }0<M<\infty\:.
\]
For each integer $M\geq1$ set 
\[
E_{M}=\{(x,h)\in X\::\: F(x,h)>M\}\:\mbox{ and }\: E_{0,M}=\cap_{N=1}^{\infty}\cup_{j=N}^{\infty}T^{-j}(E_{M}),
\]
then $\mu\times\xi_{H}(E_{M})>0$, and so $\mu\times\xi_{H}(X\setminus E_{0,M})=0$
since $\mu\times\xi_{H}$ is ergodic (see Theorem 1.5 in \cite{key-11}).
Set $\tilde{E}=D_{0}\cap(\cap_{M=1}^{\infty}E_{0,M})$, then $\mu\times\xi_{H}(X\setminus\tilde{E})=0$.
For $\xi_{H}$-a.e. $h\in H$ it holds that $\mu\{x\in K\::\:(x,h)\notin\tilde{E}\}=0$,
fix such $h_{0}\in H$ and set 
\[
A=\{x\in K\::\:(x,h_{0})\in\tilde{E}\}\:.
\]
Note that since $(x,h_{0})\in D_{0}$ for some $x\in K$, it follows
that $P_{h_{0}}\mu\ll\mathcal{H}^{m}$. It will now be shown that
\begin{equation}
\Theta^{*s-m}(\mu_{V^{h_{0}},x},x)=\infty\mbox{ for each }x\in A\:.\label{E4}
\end{equation}
Let $x\in A$, $M\geq1$ and $N\geq1$ be given, then there exists
$k\geq N$ with $T^{k}(x,h_{0})\in D_{0}\cap E_{M}$, and so $F(T^{k}(x,h_{0}))>M$.
Set $u=w_{k}(x)$ and $\beta=\left(F(x,h_{0})\right)^{-1}$, then
from Lemma \ref{L7}
\[
\mu_{V^{h_{0}},x}(K_{u})=\beta\cdot r_{u}^{s-m}\cdot F(T^{k}(x,h_{0}))\geq\beta\cdot r_{u}^{s-m}\cdot M\:.
\]
Set $d=\sup\{|y_{1}-y_{2}|\::\: y_{1},y_{2}\in K\}$, then
\[
\frac{\mu_{V^{h_{0}},x}(B(x,d\cdot r_{w_{k}(x)}))}{(2d\cdot r_{w_{k}(x)})^{s-m}}\geq\frac{\mu_{V^{h_{0}},x}(K_{u})}{(2d\cdot r_{u})^{s-m}}\geq\frac{\beta\cdot r_{u}^{s-m}\cdot M}{(2d\cdot r_{u})^{s-m}}=\frac{M\beta}{(2d)^{s-m}}\:.
\]
Since $\underset{k\rightarrow\infty}{\lim}r_{w_{k}(x)}=0$ we get
that $\Theta^{*s-m}(\mu_{V^{h_{0}},x},x)\geq\frac{M\beta}{(2d)^{s-m}}$,
and so (\ref{E4}) follows since $M$ can be chosen arbitrarily large.$\newline$Let
$x\in A$ and $y\in A\cap V_{x}^{h_{0}}$, then from (\ref{E4}) we
get
\[
\Theta^{*s-m}(\mu_{V^{h_{0}},x},y)=\Theta^{*s-m}(\mu_{V^{h_{0}},y},y)=\infty\:.
\]
Now from Theorem 6.9 in \cite{key-1} it follows that for each $M\ge1$
\[
\mathcal{H}^{s-m}(A\cap V_{x}^{h_{0}})\le M^{-1}\cdot\mu_{V^{h_{0}},x}(A\cap V_{x}^{h_{0}})\le M^{-1},
\]
and so $\mathcal{H}^{s-m}(A\cap V_{x}^{h_{0}})=0$ since $M$ can
be chosen arbitrarily large. Also, from $\mu(K\setminus A)=0$ and
Theorem 7.7 in \cite{key-1} we get that
\begin{multline*}
\int_{(V^{h_{0}})^{\perp}}\mathcal{H}^{s-m}((K\setminus A)\cap V_{y}^{h_{0}})\: d\mathcal{H}^{m}(y)\leq\\
\le const\cdot\mathcal{H}^{s}(K\setminus A)=const\cdot\mu(K\setminus A)=0\:.
\end{multline*}
This shows that $\mathcal{H}^{s-m}((K\setminus A)\cap V_{y}^{h_{0}})=0$
for $\mathcal{H}^{m}$-a.e. $y\in(V^{h_{0}})^{\perp}$, and so $\mathcal{H}^{s-m}((K\setminus A)\cap V_{x}^{h_{0}})=0$
for $\mu$-a.e. $x\in K$ since $P_{h_{0}}\mu\ll\mathcal{H}^{m}$.
It follows that for $\mu$-a.e. $x\in A$ (and so for $\mu$-a.e.
$x\in K$) we have
\[
\mathcal{H}^{s-m}(K\cap V_{x}^{h_{0}})=\mathcal{H}^{s-m}(A\cap V_{x}^{h_{0}})+\mathcal{H}^{s-m}((K\setminus A)\cap V_{x}^{h_{0}})=0\:.
\]
From this, Lemma (\ref{L4}) and Fubini's theorem, it follows that
$\mathcal{H}^{s-m}(K\cap V_{x}^{h})=0$ for $\mu\times\xi_{H}$-a.e.
$(x,h)\in K\times H$, which proves \emph{(ii)}.\emph{$\newline$$\newline$Proof
of part (iii):} Assume that $\left\Vert F_{V}\right\Vert _{\infty}<\infty$
for $\xi_{G}$-a.e. $V\in G$. From Lemma \ref{L1} and part \emph{(i)},
it follows that for $\xi_{G}$-a.e. $V\in G$ it holds for $\xi_{H}$-a.e.
$h\in H$ that 
\[
\mathcal{H}^{s-m}(K\cap(x+hV))>0\:\mbox{ for \ensuremath{\mu}-a.e. }x\in K\:.
\]
Set 
\[
B=\{(x,V)\in K\times G\::\:\mathcal{H}^{s-m}(K\cap V_{x})=0\},
\]
then from Lemma \ref{L4} we get that $B$ is a Borel set (hence universally
measurable), and so $\mu\times\xi_{G}(B)=0$ by Lemma \ref{L6}.$\newline$For
the other direction, set $\mathcal{W}=\{V\in G\::\:\left\Vert F_{V}\right\Vert _{\infty}=\infty\}$
and assume that $\xi_{G}(\mathcal{W})>0$. From part \emph{(ii)} it
follows that for $\xi_{G}$-a.e. $V\in\mathcal{W}$ we have 
\[
\mathcal{H}^{s-m}(K\cap(x+hV))=0\:\mbox{ for \ensuremath{\mu\times\xi_{H}}-a.e. }(x,h)\in X,
\]
and so from Lemma \ref{L2}
\begin{multline*}
0<\xi_{G}(\mathcal{W})\le\int\mu\times\xi_{H}\{(x,h)\::\:\mathcal{H}^{s-m}(K\cap(x+hV))=0\}\: d\xi_{G}(V)=\\
=\int\int\xi_{H}\{h\::\:\mathcal{H}^{s-m}(K\cap(x+hgU))=0\}\: d\xi_{O}(g)\: d\mu(x)=\\
=\int\int\xi_{H}(\{h\::\:\mathcal{H}^{s-m}(K\cap(x+hU))=0\}\cdot g^{-1})\: d\xi_{O}(g)\: d\mu(x)=\\
=\int\xi_{O}\{g\::\:\mathcal{H}^{s-m}(K\cap(x+gU))=0\}\: d\mu(x)=\\
=\int\xi_{G}\{V\::\:\mathcal{H}^{s-m}(K\cap V_{x})=0\}\: d\mu(x)=\\
=\mu\times\xi_{G}\{(x,V)\::\:\mathcal{H}^{s-m}(K\cap V_{x})=0\},
\end{multline*}
which completes the proof of \emph{(iii)}.$\newline$Part \emph{(iv)}
can be proven in a similar manner, and so the proof of Theorem \ref{T4}
is complete. $\square$

\emph{Proof of theorem \ref{T2}, part (i)}: Let $M>0$ be so large
such that for 
\[
E=\{(x,h)\in X\::\: F(x,h)\leq M\}
\]
we have $\mu\times\xi_{H}(E)>0$. Set $E_{0}=\cap_{N=1}^{\infty}\cup_{j=N}^{\infty}T^{-j}(E)$,
then $\mu\times\xi_{H}(X\setminus E_{0})=0$ since $\mu\times\xi_{H}$
is ergodic. Set $E_{1}=E_{0}\cap D_{0}$, then $\mu\times\xi_{H}(X\setminus E_{1})=0$.
For $\xi_{H}$-a.e. $h\in H$ it holds that $\mu\{x\in K\::\:(x,h)\notin E_{1}\}=0$,
fix such $h_{0}\in H$. For each $l\geq1$ set
\[
A_{l}=\{x\in K\::\:(x,h_{0})\in E_{1}\mbox{ and }F(x,h_{0})\geq l^{-1}\},
\]
and fix $l_{0}\geq1$. It will now be shown that
\begin{equation}
\Theta_{*}^{s-m}(\mu_{V^{h_{0}},x},x)\leq(2\rho)^{m-s}l_{0}M\:\mbox{ for each }x\in A_{l_{0}}\:.\label{E1}
\end{equation}
Let $x\in A_{l_{0}}$ and let $N\geq1$ be given, then since $(x,h_{0})\in E_{1}$
it follows that there exist $k\geq N$ with $T^{k}(x,h_{0})\in E\cap D_{0}$,
and so $F(T^{k}(x,h_{0}))\leq M$. Set $u=w_{k}(x)$, then from Lemma
\ref{L7} we have
\[
\mu_{V^{h_{0}},x}(K_{u})=\left(F(x,h_{0})\right)^{-1}\cdot r_{u}^{s-m}\cdot F(T^{k}(x,h_{0}))\leq l_{0}r_{u}^{s-m}M\:,
\]
from which it follows that
\[
\frac{\mu_{V^{h_{0}},x}(B(x,\rho\cdot r_{w_{k}(x)}))}{(2\rho\cdot r_{w_{k}(x)})^{s-m}}\leq\frac{\mu_{V^{h_{0}},x}(K_{u})}{(2\rho\cdot r_{u})^{s-m}}\leq\frac{l_{0}r_{u}^{s-m}M}{(2\rho\cdot r_{u})^{s-m}}=(2\rho)^{m-s}l_{0}M\:.
\]
This proves (\ref{E1}) since $r_{w_{k}(x)}$ tends to $0$ as $k$
tends to $\infty$.$\newline$$\newline$As in the proof of part \emph{(i)}
of Theorem \ref{T4}, from $\mu(K\setminus\cup_{l=1}^{\infty}A_{l})=0$
it follows that for $\mu$-a.e. $x\in K$ there exists $l_{x}\geq1$
with $\mu_{V^{h_{0}},x}(A_{l_{x}}\cap V_{x}^{h_{0}})>0$. Fix such
an $x_{0}$ and let $y\in A_{l_{x_{0}}}\cap V_{x_{0}}^{h_{0}}$, then
from (\ref{E1}) we get 
\[
\Theta_{*}^{s-m}(\mu_{V^{h_{0}},x_{0}},y)=\Theta_{*}^{s-m}(\mu_{V^{h_{0}},y},y)\leq(2\rho)^{m-s}l_{x_{0}}M,
\]
and so from Theorem 6.11 in \cite{key-1} it follows that 
\[
\mathcal{P}^{s-m}(K\cap V_{x_{0}}^{h_{0}})\geq\mathcal{P}^{s-m}(A_{l_{x_{0}}}\cap V_{x_{0}}^{h_{0}})\geq(2\rho)^{s-m}l_{x_{0}}^{-1}M^{-1}\cdot\mu_{V^{h_{0}},x_{0}}(A_{l_{x_{0}}}\cap V_{x_{0}}^{h_{0}})>0.
\]
Since $\xi_{G}(G\setminus\mathcal{V})=0$, this shows that for $\xi_{G}$-a.e.
$V\in G$ it holds for $\xi_{H}$-a.e. $h\in H$ that $\mathcal{P}^{s-m}(K\cap(x+hV))>0$
for $\mu$-a.e. $x\in K$. Set 
\[
B=\{(x,V)\in K\times G\::\:\mathcal{P}^{s-m}(K\cap V_{x})=0\},
\]
then from Lemma \ref{L5} we get that $B$ is universally measurable,
and so the claim stated in \emph{(i)} follows from Lemma \ref{L6}.$\newline$$\newline$\emph{Proof
of part (ii):} Assume $V$ is such that $\left\Vert \frac{1}{F}\right\Vert _{L^{\infty}(\mu\times\xi_{H})}=\infty$,
then 
\[
\mu\times\xi_{H}\{(x,h)\::\: F(x,h)<M^{-1}\}>0\:\mbox{ for each }0<M<\infty\:.
\]
For each integer $M\geq1$ set 
\[
E_{M}=\{(x,h)\::\: F(x,h)<M^{-1}\}\mbox{ and }E_{0,M}=\cap_{N=1}^{\infty}\cup_{j=N}^{\infty}T^{-j}(E_{M}),
\]
then since $\mu\times\xi_{H}$ is ergodic and $\mu\times\xi_{H}(E_{M})>0$
it follows that $\mu\times\xi_{H}(X\setminus E_{0,M})=0$. Set $\tilde{E}=D_{0}\cap(\cap_{M=1}^{\infty}E_{0,M})$,
then $\mu\times\xi_{H}(X\setminus\tilde{E})=0$. For $\xi_{H}$-a.e.
$h\in H$ it holds that $\mu\{x\in K\::\:(x,h)\notin\tilde{E}\}=0$,
fix such $h_{0}\in H$ and set $A=\{x\in K\::\:(x,h_{0})\in\tilde{E})\}$.$\newline$It
will now be shown that 
\begin{equation}
\Theta_{*}^{s-m}(\mu_{V^{h_{0}},x},x)=0\mbox{ for each }x\in A\:.\label{E3}
\end{equation}
Let $x\in A$, $M\geq1$ and $N\geq1$ be given, then there exists
$k\geq N$ with $T^{k}(x,h_{0})\in D_{0}\cap E_{M}$, and so $F(T^{k}(x,h_{0}))<M^{-1}$.
Set $u=w_{k}(x)$, then from Lemma \ref{L7}
\[
\mu_{V^{h_{0}},x}(K_{u})=\left(F(x,h_{0})\right)^{-1}\cdot r_{u}^{s-m}\cdot F(T^{k}(x,h_{0}))\leq\left(F(x,h_{0})\right)^{-1}\cdot r_{u}^{s-m}\cdot M^{-1}\:,
\]
from which it follows that
\begin{multline*}
\frac{\mu_{V^{h_{0}},x}(B(x,\rho\cdot r_{w_{k}(x)}))}{(2\rho\cdot r_{w_{k}(x)})^{s-m}}\leq\frac{\mu_{V^{h_{0}},x}(K_{u})}{(2\rho\cdot r_{u})^{s-m}}\leq\\
\leq\frac{\left(F(x,h_{0})\right)^{-1}\cdot r_{u}^{s-m}\cdot M^{-1}}{(2\rho\cdot r_{u})^{s-m}}=(2\rho)^{m-s}\cdot\left(F(x,h_{0})\right)^{-1}\cdot M^{-1}.
\end{multline*}
This shows that 
\[
\Theta_{*}^{s-m}(\mu_{V^{h_{0}},x},x)\leq(2\rho)^{m-s}\cdot\left(F(x,h_{0})\right)^{-1}\cdot M^{-1},
\]
and so (\ref{E3}) holds since $M$ can be chosen arbitrarily large.$\newline$We
have
\[
0=\mu(K\setminus A)=\int\mu_{V^{h_{0}},x}(K\setminus A)\: d\mu(x),
\]
hence $\mu_{V^{h_{0}},x}(A\cap V_{x}^{h_{0}})>0$ for $\mu$-a.e.
$x\in K$. Fix such $x_{0}\in K$ and let $y\in A\cap V_{x_{0}}^{h_{0}}$,
then from (\ref{E3}) we get 
\[
\Theta_{*}^{s-m}(\mu_{V^{h_{0}},x_{0}},y)=\Theta_{*}^{s-m}(\mu_{V^{h_{0}},y},y)=0\:.
\]
Now from Theorem 6.11 in \cite{key-1} it follows that for each $\epsilon>0$
\[
\mathcal{P}^{s-m}(K\cap V_{x_{0}}^{h_{0}})\geq\mathcal{P}^{s-m}(A\cap V_{x_{0}}^{h_{0}})\geq\epsilon^{-1}\cdot\mu_{V^{h_{0}},x_{0}}(A\cap V_{x_{0}}^{h_{0}}),
\]
which shows that $\mathcal{P}^{s-m}(K\cap V_{x_{0}}^{h_{0}})=\infty$
since $\epsilon$ can be chosen arbitrarily small and $\mu_{V^{h_{0}},x_{0}}(A\cap V_{x_{0}}^{h_{0}})>0$.$\newline$This
proves that if $\left\Vert \frac{1}{F_{V}}\right\Vert _{L^{\infty}(\mu\times\xi_{H})}=\infty$,
then for $\xi_{H}$-a.e. $h\in H$ we have $\mathcal{P}^{s-m}(K\cap(x+hV))=\infty$
for $\mu$-a.e. $x\in K$, and so \emph{(ii)} follows from Lemma (\ref{L5})
and Fubini's theorem.

\emph{Proof of part (iii):} Assume that $\left\Vert \frac{1}{F_{V}}\right\Vert _{L^{\infty}(\mu\times\xi_{H})}=\infty$
for $\xi_{G}$-a.e. $V\in G$, then from Lemma \ref{L1} and part
\emph{(ii)} it follows that for $\xi_{G}$-a.e. $V\in G$ it holds
for $\xi_{H}$-a.e. $h\in H$ that $\mathcal{P}^{s-m}(K\cap(x+hV))=\infty$
for $\mu$-a.e. $x\in K$. Set 
\[
B=\{(x,V)\in K\times G\::\:\mathcal{P}^{s-m}(K\cap V_{x})<\infty\}\:,
\]
then from Lemma \ref{L5} we get that $B$ is universally measurable,
and so the claim stated in \emph{(iii)} follows from Lemma \ref{L6}.
This completes the proof of Theorem \ref{T2}. $\square$

\subsection{Proofs of Corollaries \ref{C2}, \ref{C3} and \ref{C1}}

The following lemma will be used in the proofs of Corollaries \ref{C2}
and \ref{C1}. For its proof see Lemma 3.2 in \cite{key-9} and the
discussion before it.
\begin{lem}
\label{L12}Assume $m=1$ and $s>2$, then $P_{V^{\perp}}\mu\ll\mathcal{H}^{m}$
and $\frac{dP_{V^{\perp}}\mu}{d\mathcal{H}^{m}}$ has a continuous
version for $\xi_{G}$-a.e. $V\in G$.
\end{lem}

\emph{Proof of corollary }\ref{C2}\emph{:} Assuming $m=1$, $s>2$
and $|H|<\infty$, it will be shown that $\left\Vert F_{V}\right\Vert _{L^{\infty}(\mu\times\xi_{H})}<\infty$
for $\xi_{G}$-a.e. $V\in G$. From this and from part \emph{(iii)}
of Theorem \ref{T4} the corollary will follow. Set
\[
E=\{V\in G\::\: P_{V^{\perp}}\mu\ll\mathcal{H}^{m}\mbox{ and }\frac{dP_{V^{\perp}}\mu}{d\mathcal{H}^{m}}\mbox{ is continuous}\},
\]
then from Lemma \ref{L12} we get $\xi_{G}(G\setminus E)=0$. From
this and from Lemma \ref{L2} it now follows that
\begin{multline*}
0=\xi_{G}(G\setminus E)=\xi_{O}\{g\in O(n)\::\: gU\notin E\}=\\
=\int\xi_{H}\{h\::\: hgU\notin E\}\: d\xi_{O}(g)=\int\xi_{H}\{h\::\: hV\notin E\}\: d\xi_{G}(V),
\end{multline*}
and so $\xi_{H}\{h\::\: hV\notin E\}=0$ for $\xi_{G}$-a.e. $V$.
We fix such a $V\in G$. Since $|H|<\infty$, for each $h\in H$ we
have $\xi_{H}\{h\}>0$, and so $hV\in E$.$\newline$For each $h\in H$
and $y\in(hV)^{\perp}$ set $Q_{h}(y)=\Theta_{*}^{m}(P_{(hV)^{\perp}}\mu,y)$,
fix $h_{0}\in H$, and set $W=(h_{0}V)^{\perp}$. Since $\mathcal{H}^{m}(B(y,r)\cap W)=(2\epsilon)^{m}$
for each $y\in W$ and $0<\epsilon<\infty$, it follows from Theorem
2.12 in \cite{key-1} that $Q_{h_{0}}(y)=\frac{dP_{W}\mu}{d\mathcal{H}^{m}}(y)$
for $\mathcal{H}^{m}$-a.e. $y\in W$, i.e. the function $Q_{h_{0}}$
equals a continuous function as members of $L^{1}(W,\mathcal{H}^{m})$.
Also, since $\mu$ is supported on a compact set it follows that the
set $\{y\in W\::\: Q_{h_{0}}(y)\ne0\}$ is bounded, so $Q_{h_{0}}$
equals a continuous function with compact support in $L^{1}(W,\mathcal{H}^{m})$,
which shows that $\left\Vert Q_{h_{0}}\right\Vert _{L^{\infty}(W,\mathcal{H}^{m})}<\infty$.
Since $P_{W}\mu\ll\mathcal{H}^{m}$ it follows that $\left\Vert Q_{h_{0}}\right\Vert _{L^{\infty}(P_{W}\mu)}<\infty$.$\newline$Now
set $M=\max\{\left\Vert Q_{h}\right\Vert _{L^{\infty}(P_{(hV)^{\perp}}\mu)}\::\: h\in H\}$,
then $M<\infty$ since $|H|<\infty$. Also, we have
\begin{multline*}
0=\frac{1}{|H|}\sum_{h\in H}P_{(hV)^{\perp}}\mu\{y\in(hV)^{\perp}\::\:|Q_{h}(y)|>M\}=\\
=\frac{1}{|H|}\sum_{h\in H}\mu\{x\in K\::\:|Q_{h}(P_{(hV)^{\perp}}(x))|>M\}=\\
=\frac{1}{|H|}\sum_{h\in H}\mu\{x\in K\::\:|F_{V}(x,h)|>M\}=\\
=\int\mu\{x\in K\::\:|F_{V}(x,h)|>M\}\: d\xi_{H}(h)=\\
=\mu\times\xi_{H}\{(x,h)\in K\times H\::\:|F_{V}(x,h)|>M\},
\end{multline*}
which shows that $\left\Vert F_{V}\right\Vert _{L^{\infty}(\mu\times\xi_{H})}\leq M<\infty$.
This completes the proof of corollary \ref{C2}. $\square$\emph{$\newline$$\newline$Proof
of corollary }\ref{C3}\emph{:} Assume that $H=O(n)$ and
\[
\mu\times\xi_{G}\{(x,V)\::\:\mathcal{H}^{s-m}(K\cap V_{x})>0\}>0\:.
\]
Let $V\in\mathcal{V}$, then since $\xi_{H}=\xi_{O}$ we have
\[
\mu\times\xi_{H}\{(x,h)\::\:\mathcal{H}^{s-m}(K\cap(x+hV))>0\}>0\:,
\]
and so from part \emph{(ii)} of theorem \ref{T4} it follows that
$\left\Vert F_{V}\right\Vert _{L^{\infty}(\mu\times\xi_{H})}<\infty$.
Set $M=\left\Vert F_{V}\right\Vert _{L^{\infty}(\mu\times\xi_{H})}$,
set
\[
E=\{W\in G\::\: P_{W^{\perp}}\mu\ll\mathcal{H}^{m}\mbox{ and }\left\Vert \frac{dP_{W^{\perp}}\mu}{d\mathcal{H}^{m}}\right\Vert _{L^{\infty}(\mathcal{H}^{m})}\le M\},
\]
and for each $h\in H$ set $P_{h}=P_{(hV)^{\perp}}$.$\newline$We
shall first show that $\xi_{G}(G\setminus E)=0$. Since $P_{W^{\perp}}\mu\ll\mathcal{H}^{m}$
for $\xi_{G}$-a.e. $W\in G$ (see the proof of lemma \ref{L1}),
and since $\xi_{H}=\xi_{O}$, we have 
\begin{multline}
\xi_{G}(G\setminus E)=\xi_{G}(G\setminus\{W\in G\::\: P_{W^{\perp}}\mu\ll\mathcal{H}^{m}\})+\\
+\xi_{G}\{W\in G\::\: P_{W^{\perp}}\mu\ll\mathcal{H}^{m}\mbox{ and }\left\Vert \frac{dP_{W^{\perp}}\mu}{d\mathcal{H}^{m}}\right\Vert _{L^{\infty}(\mathcal{H}^{m})}>M\}=\\
=\xi_{H}\{h\::\: P_{h}\mu\ll\mathcal{H}^{m}\mbox{ and }\left\Vert \frac{dP_{h}\mu}{d\mathcal{H}^{m}}\right\Vert _{L^{\infty}(\mathcal{H}^{m})}>M\}\:.\label{E16}
\end{multline}
Let $h\in H$ be such that $P_{h}\mu\ll\mathcal{H}^{m}$ and $\left\Vert \frac{dP_{h}\mu}{d\mathcal{H}^{m}}\right\Vert _{L^{\infty}(P_{h}\mu)}\le M$,
then
\begin{multline*}
0=P_{h}\mu\{y\in(hV)^{\perp}\::\:\frac{dP_{h}\mu}{d\mathcal{H}^{m}}(y)>M\}=\\
=\int_{(hV)^{\perp}}1_{\{\frac{dP_{h}\mu}{d\mathcal{H}^{m}}>M\}}\cdot\frac{dP_{h}\mu}{d\mathcal{H}^{m}}\: d\mathcal{H}^{m}\ge\\
\ge M\cdot\mathcal{H}^{m}\{y\in(hV)^{\perp}\::\:\frac{dP_{h}\mu}{d\mathcal{H}^{m}}(y)>M\},
\end{multline*}
which shows that $\left\Vert \frac{dP_{h}\mu}{d\mathcal{H}^{m}}\right\Vert _{L^{\infty}(\mathcal{H}^{m})}\le M$.
From this and from (\ref{E16}) it follows that
\begin{equation}
\xi_{G}(G\setminus E)=\xi_{H}\{h\::\: P_{h}\mu\ll\mathcal{H}^{m}\mbox{ and }\left\Vert \frac{dP_{h}\mu}{d\mathcal{H}^{m}}\right\Vert _{L^{\infty}(P_{h}\mu)}>M\}\:.\label{E17}
\end{equation}
From Theorem 2.12 in \cite{key-1} we get that for each $h\in H$
with $P_{h}\mu\ll\mathcal{H}^{m}$ 
\[
F_{V}(x,h)=\frac{dP_{h}\mu}{d\mathcal{H}^{m}}(P_{h}(x))\:\mbox{ for \ensuremath{\mu}-a.e. \ensuremath{x\in K}},
\]
and so from (\ref{E17}) 
\begin{multline*}
\xi_{G}(G\setminus E)\le\xi_{H}\{h\::\:\left\Vert F_{V}(\cdot,h)\right\Vert _{L^{\infty}(\mu)}>M\}=\\
=\xi_{H}\{h\::\:\mu\{x\::\: F_{V}(x,h)>\left\Vert F_{V}\right\Vert _{L^{\infty}(\mu\times\xi_{H})}\}>0\}=0\:.
\end{multline*}
Since $\xi_{G}(\mathcal{W})>0$ for every non-empty open set $\mathcal{W}\subset G$,
it follows from $\xi_{G}(G\setminus E)=0$ that $E$ is dense in $G$,
and so in order to prove the corollary it suffice to show that $E$
is a closed subset of $G$. Let $W_{0}\in\overline{E}$, let $y\in W_{0}^{\perp}$
and let $r\in(0,\infty)$. Given $\epsilon>0$ there exists $W\in E$
so close to $W_{0}$ in $G$ (with respect to the metric $d_{G}$
defined in section \ref{S3.1}), such that
\[
P_{W_{0}^{\perp}}^{-1}(B(y,r))\cap K\subset P_{W^{\perp}}^{-1}(B(P_{W^{\perp}}y,r+\epsilon)).
\]
From this and since $W\in E$ it follows that
\begin{multline*}
P_{W_{0}^{\perp}}\mu(B(y,r))=\mu(P_{W_{0}^{\perp}}^{-1}(B(y,r))\cap K)\le\\
\le\mu(P_{W^{\perp}}^{-1}(B(P_{W^{\perp}}y,r+\epsilon)))=P_{W^{\perp}}\mu((B(P_{W^{\perp}}y,r+\epsilon)))=\\
=\int_{B(P_{W^{\perp}}y,r+\epsilon)\cap W^{\perp}}\frac{dP_{W^{\perp}}\mu}{d\mathcal{H}^{m}}\: d\mathcal{H}^{m}\le\\
\le M\cdot\mathcal{H}^{m}(B(P_{W^{\perp}}y,r+\epsilon)\cap W^{\perp})=M\cdot(2\cdot(r+\epsilon))^{m},
\end{multline*}
and since this holds for each $\epsilon>0$ we have
\[
P_{W_{0}^{\perp}}\mu(B(y,r)\cap W_{0}^{\perp})\le M\cdot(2r)^{m}=M\cdot\mathcal{H}^{m}(B(y,r)\cap W_{0}^{\perp})\:.
\]
This holds for every $y\in W_{0}^{\perp}$ and $r\in(0,\infty)$,
hence $W_{0}\in E$ by Theorem 2.12 in \cite{key-1}, which shows
that $E$ is closed in $G$ and completes the proof of the corollary.
$\square$$\newline$$\newline$\emph{Proof of corollary }\ref{C1}\emph{:}
Assuming $m=1$ and $s>2$, it will be shown that $\left\Vert \frac{1}{F_{V}}\right\Vert _{L^{\infty}(\mu\times\xi_{H})}=\infty$
for $\xi_{G}$-a.e. $V\in G$. From this and part \emph{(iii)} of
Theorem \ref{T2} the corollary will follow. Set
\[
E=\{V\in G\::\: P_{V^{\perp}}\mu\ll\mathcal{H}^{m}\mbox{ and }\frac{dP_{V^{\perp}}\mu}{d\mathcal{H}^{m}}\mbox{ is continuous}\},
\]
then as in the proof of corollary \ref{C2} it follows from Lemma
\ref{L12} and Lemma \ref{L2} that
\[
0=\xi_{G}(G\setminus E)=\int\xi_{H}\{h\::\: hV\notin E\}\: d\xi_{G}(V),
\]
and so $\xi_{H}\{h\::\: hV\notin E\}=0$ for $\xi_{G}$-a.e. $V$.
Fix such $V\in G$, let $M>0$, set 
\[
A=\{h\in H\::\: hV\in E\},
\]
and for each $h\in H$ and $y\in(hV)^{\perp}$ set $Q_{h}(y)=\Theta_{*}^{m}(P_{(hV)^{\perp}}\mu,y)$
and 
\[
L_{h}=\{y\in(hV)^{\perp}\::\:0<Q_{h}(y)\leq M^{-1}\}.
\]
Fix $h_{0}\in A$ and set $W=(h_{0}V)^{\perp}$. From Theorem 2.12
in \cite{key-1} it follows that $Q_{h_{0}}(y)=\frac{dP_{W}\mu}{d\mathcal{H}^{m}}(y)$
for $\mathcal{H}^{m}$-a.e. $y\in W$, hence the function $Q_{h_{0}}$
equals a continuous function in $L^{1}(W,\mathcal{H}^{m})$. Also,
since $\mu$ is supported on a compact set, it follows that the set
$\{y\in W\::\: Q_{h_{0}}(y)\ne0\}$ is bounded. From these two facts
it easily follows that $\mathcal{H}^{m}(L_{h_{0}})>0$, and so $P_{W}\mu(L_{h_{0}})>0$
since $Q_{h_{0}}=\frac{dP_{W}\mu}{d\mathcal{H}^{m}}$ and $Q_{h_{0}}>0$
on $L_{h_{0}}$. From this we get that
\[
0<\mu\{x\in K\::\: Q_{h_{0}}(P_{W}(x))\leq M^{-1}\}=\mu\{x\in K\::\: F_{V}(x,h_{0})\leq M^{-1}\},
\]
and so by Fubini's theorem 
\[
\mu\times\xi_{H}\{(x,h)\::\:\frac{1}{F_{V}(x,h)}\ge M\}=\int_{A}\mu\{x\in K\::\: F_{V}(x,h)\leq M^{-1}\}\: d\xi_{H}(h)>0.
\]
It follows that $\left\Vert \frac{1}{F_{V}}\right\Vert _{L^{\infty}(\mu\times\xi_{H})}\geq M$,
and so $\left\Vert \frac{1}{F_{V}}\right\Vert _{L^{\infty}(\mu\times\xi_{H})}=\infty$
since we can choose $M$ as large as we want. This completes the proof
of the corollary. $\square$$\newline$

\section{\label{S4}Proof of Theorem \ref{T6}}

Set $\Lambda=\{1,2\}$. Given $0<\rho<\frac{1}{2}$, define $f_{\rho,1},f_{\rho,2}:\mathbb{R}\rightarrow\mathbb{R}$
by $f_{\rho,1}(x)=\rho\cdot x$ and $f_{\rho,2}(x)=\rho\cdot x+1-\rho$
for each $x\in\mathbb{R}$, let $C_{\rho}\subset[0,1]$ be the attractor
of the IFS $\{f_{\rho,1},f_{\rho,2}\}$, set $d_{\rho}=\dim_{H}C_{\rho}$
(so that $d_{\rho}=\frac{\log2}{\log\rho^{-1}}$), and for each $E\subset\mathbb{R}$
set $\mu_{\rho}(E)=\frac{\mathcal{H}^{d_{\rho}}(C_{\rho}\cap E)}{\mathcal{H}^{d_{\rho}}(C_{\rho})}$.
Let $0<a<b<\frac{1}{2}$ be such that $\frac{1}{a}$ and $\frac{1}{b}$
are Pisot numbers, $\frac{\log b}{\log a}$ is irrational, and $d_{a}+d_{b}>1$.
Let $I=[0,1)$ and let $\mathcal{L}$ be Lebesgue measure on $I$.
Fix $\tau\in(0,\infty)$, and for each $t\in I$ and $z\in\mathbb{R}^{2}$
define $W^{t}=\{x\cdot(1,\tau\cdot a^{t})\::\: x\in\mathbb{R}\}$,
$V^{t}=(W^{t})^{\perp}$ and $V_{z}^{t}=z+V^{t}$. In order to prove
Theorem \ref{T6} we shall first prove the following:
\begin{thm}
\label{T13}For $\mu_{a}\times\mu_{b}\times\mathcal{L}$-a.e. $(x,y,t)\in C_{a}\times C_{b}\times I$
it holds that 
\[
\mathcal{H}^{d_{a}+d_{b}-1}((C_{a}\times C_{b})\cap V_{(x,y)}^{t})=0\:.
\]

\end{thm}

\subsection{Preliminaries}

Set $\alpha=\frac{\log b}{\log a}$ (so $\alpha\in I\setminus\mathbb{Q}$),
and for each $t\in I$ set $R(t)=t+\alpha\mod1$. Given $0<\rho<\frac{1}{2}$
and a word $\lambda_{1}\cdot...\cdot\lambda_{l}=w\in\Lambda^{*}$,
write $f_{\rho,w}=f_{\rho,\lambda_{1}}\circ...\circ f_{\rho,\lambda_{l}}$
and $C_{\rho,w}=f_{\rho,w}(C_{\rho})$. For each $n\geq1$ and $x\in C_{\rho}$
let $w_{\rho,n}(x)\in\Lambda^{n}$ be the unique word of length $n$
which satisfies $x\in C_{\rho,w_{\rho,n}(x)}$, and let $S_{\rho}(x)=f_{\rho,w_{\rho,1}(x)}^{-1}(x)$.
We also write $w_{\rho,0}(x)=\emptyset$ and $C_{\rho,\emptyset}=C_{\rho}$. 

The following dynamical system will be used in the proof of Theorem
\ref{T13}. The idea of using this system comes from the partition
operator introduced in section 10 of \cite{key-10}. Set $K=C_{a}\times C_{b}$,
$X=K\times I$, $\mu=\mu_{a}\times\mu_{b}$, $\nu=\mu\times\mathcal{L}$,
and for each $(x,y,t)\in X$ define
\[
T(x,y,t)=\begin{cases}
(x,S_{b}(y),R(t)) & \mbox{, if }t\in[0,1-\alpha)\\
((S_{a}(x),S_{b}(y),R(t)) & \mbox{, else}
\end{cases}\:.
\]
It is easy to check that the system $(X,\nu,T)$ is measure preserving,
and from Lemma 2.2 in \cite{key-16} it follows that it is ergodic.

Let $\mathcal{R}$ be the Borel $\sigma$-algebra of $\mathbb{R}^{2}$.
For each $t\in I$ let $P_{t}$ be the orthogonal projection onto
$W^{t}$, and let $\{\mu_{t,z}\}_{z\in\mathbb{R}^{2}}$ be the disintegration
of $\mu$ with respect to $P_{t}^{-1}(\mathcal{R})$ (see section
\ref{S3.1} above). Also, for each $(z,t)\in X$ define $F(z,t)=\Theta_{*}^{1}(P_{t}\mu,P_{t}z)$.

\subsection{Auxiliary lemmas}
\begin{lem}
\label{L8}It holds that $I_{1}(\mu)<\infty$, where $I_{1}(\mu)$
is the $1$-energy of $\mu$.
\end{lem}

\emph{Proof:} Set $\delta=1-2b$, then for each $(x,y)\in\mathbb{R}^{2}$
and $k\geq1$
\begin{multline*}
\mu(B((x,y),\delta\cdot a^{k}))\leq\mu((x-\delta\cdot a^{k},x+\delta\cdot a^{k})\times(y-\delta\cdot a^{k},y+\delta\cdot a^{k}))\leq\\
\le\mu_{a}(x-\delta\cdot a^{k},x+\delta\cdot a^{k})\cdot\mu_{b}(y-\delta\cdot a^{k},y+\delta\cdot a^{k})\le2^{-k}\cdot2^{-[k\log_{b}a]}\le\\
\le2^{-k}\cdot2^{1-k\log_{b}a}=2\cdot a^{k(1+\log_{b}a)\log_{a}2^{-1}}=2\cdot a^{k(d_{a}+d_{b})}\:.
\end{multline*}
This shows that there exists a constant $M>0$ with $\mu(B(z,r))\leq M\cdot r^{d_{a}+d_{b}}$
for each $z\in\mathbb{R}^{2}$ and $r>0$. Since $d_{a}+d_{b}>1$,
the lemma follows from the discussion found at the beginning of chapter
8 of \cite{key-1}. $\square$
\begin{lem}
\label{L9}Let $n_{1},n_{2}\geq1$, $w_{1}\in\Lambda^{n_{1}}$ and
$w_{2}\in\Lambda^{n_{2}}$. For each $(x,y)\in K$ set $g(x,y)=(f_{a,w_{1}}(x),f_{b,w_{2}}(y))$,
then for each Borel set $B\subset K$
\[
\mu(g(B))=2^{-n_{1}-n_{2}}\cdot\mu(B)\:.
\]

\end{lem}

\emph{Proof:} We prove this by using the $\pi-\lambda$ theorem (see
\cite{key-15}). Let $\mathcal{E}$ be the collection of all Borel
sets $B\subset K$ which satisfy $\mu(g(B))=2^{-n_{1}-n_{2}}\cdot\mu(B)$,
then $\mathcal{E}$ is a $\lambda$-system. Set 
\[
\mathcal{P}=\{C_{a,u_{1}}\times C_{b,u_{2}}\::\: u_{1},u_{2}\in\Lambda^{*}\}\cup\{\emptyset\},
\]
then $\mathcal{P}$ is a $\pi$-system, $\mathcal{P}\subset\mathcal{E}$
and $\sigma(\mathcal{P})$ equals the collection of all Borel subsets
of $K$. From the $\pi-\lambda$ theorem it follows that $\sigma(\mathcal{P})\subset\mathcal{E}$,
hence $\mathcal{E}$ equals the collection of all Borel subsets of
$K$, and the lemma is proven. $\square$
\begin{lem}
\label{L10}It holds that $0<\mathcal{H}^{d_{a}+d_{b}}(K)<\infty$,
and $\mu(E)=\frac{\mathcal{H}^{d_{a}+d_{b}}(K\cap E)}{\mathcal{H}^{d_{a}+d_{b}}(K)}$
for each Borel set $E\subset\mathbb{R}^{2}$.
\end{lem}

\emph{Proof:} From Theorem 8.10 in \cite{key-1} it follows that $\mathcal{H}^{d_{a}+d_{b}}(K)>0$,
and by an elementary covering argument it can be shown that $\mathcal{H}^{d_{a}+d_{b}}(K)<\infty$.
The rest of the lemma can be proven by using the $\pi-\lambda$ theorem,
as in the proof of Lemma \ref{L9}. $\square$
\begin{lem}
\label{L11}Let $0<M<\infty$ and set $E_{M}=\{(z,t)\in X\::\: F(z,t)>M\}$,
then $\nu(E_{M})>0$.
\end{lem}

\emph{Proof:} Assume by contradiction that $\nu(E_{M})=0$ and set
\[
L=\{t\in I\::\:\mu\{z\::\:(z,t)\in E_{M}\}=0\},
\]
then $\mathcal{L}(I\setminus L)=0$, and so $\overline{L}=I$. Set
\[
A=\{t\in I\::\: P_{t}\mu\ll\mathcal{H}^{1}\mbox{ and }\left\Vert \frac{dP_{t}\mu}{d\mathcal{H}^{1}}\right\Vert _{L^{\infty}(\mathcal{H}^{1})}\le M\},
\]
and let $t\in L$. For $P_{t}\mu$-a.e. $z\in W^{t}$ we have $\Theta_{*}^{1}(P_{t}\mu,z)\le M$,
hence from parts (2) and (3) of Theorem 2.12 in \cite{key-1} it follows
that $t\in A$. This shows that $L\subset A$, and so that $\overline{A}=I$.
By an argument similar to the one given at the end of the proof of
Corollary \ref{C3}, it can be shown that $A$ is a closed subset
of $I$, and so $A=I$. In particular it follows that $P_{t}\mu\ll\mathcal{H}^{1}$
for each $t\in I$, which is a contradiction to Theorem 4.1 in \cite{key-7}.
This shows that we must have $\nu(E_{M})>0$, and the lemma is proven.
$\square$

\subsection{Proofs of Theorems \ref{T13} and \ref{T6}}

\emph{Proof of theorem }\ref{T13}\emph{:} Let $D$ be the set of
all $(z,t)\in X$ such that $P_{t}\mu\ll\mathcal{H}^{1}$, $\mu_{t,z}$
is defined,
\[
\mu_{t,z}(C_{a,w_{1}}\times C_{b,w_{2}})=\underset{\epsilon\downarrow0}{\lim}\frac{\mu((C_{a,w_{1}}\times C_{b,w_{2}})\cap P_{t}^{-1}(B(P_{t}z,\epsilon)))}{P_{t}\mu(B(P_{t}z,\epsilon))}
\]
for each $w_{1},w_{2}\in\Lambda^{*}$, and
\[
0<F(z,t)=\underset{\epsilon\downarrow0}{\lim}\frac{P_{t}\mu(B(P_{t}z,\epsilon))}{2\epsilon}<\infty\:.
\]
From Lemma \ref{L8} and from the same arguments as the ones given
at the beginning of section \ref{S3.3}, it follows that $\nu(X\setminus D)=0$.
Set $D_{0}=\cap_{j=0}^{\infty}T^{-j}D$, then $\nu(X\setminus D_{0})=0$
since $T$ is measure preserving.$\newline$For $0<M<\infty$ let
$E_{M}$ be as in Lemma \ref{L11}, and set $E_{0,M}=\cap_{N=1}^{\infty}\cup_{j=N}^{\infty}T^{-j}(E_{M})$.
Since $\nu(E_{M})>0$, it follows from the ergodicity of $(X,\nu,T)$
that $\nu(X\setminus E_{0,M})=0$. Set $D_{1}=D_{0}\cap(\cap_{M=1}^{\infty}E_{0,M})$,
then $\nu(X\setminus D_{1})=0$. For $\mathcal{L}$-a.e. $t\in I$
it holds that $\mu\{z\in K\::\:(z,t)\notin D_{1}\}=0$, fix such $t_{0}\in I$
and set $A=\{z\in K\::\:(z,t_{0})\in D_{1}\}$. Note that from $A\ne\emptyset$
it follows that $P_{t_{0}}\mu\ll\mathcal{H}^{1}$.$\newline$Set $\eta=d_{a}+d_{b}-1$.
It will now be shown that 
\begin{equation}
\Theta^{*\eta}(\mu_{t_{0},z},z)=\infty\;\mbox{ for each }z\in A\:.\label{E7}
\end{equation}
Let $(x,y)=z\in A$ and set $\beta=\left(F(z,t_{0})\right)^{-1}$,
then $0<\beta<\infty$ since $(z,t_{0})\in D_{0}$. Let $M\geq1$
and $N\geq1$ be given, then there exists $k\geq N$ with $T^{k}(z,t_{0})\in D_{0}\cap E_{M}$,
and so $F(T^{k}(z,t_{0}))>M$. Set $l=[t_{0}+k\alpha]$, then
\begin{multline}
\mu_{t_{0},z}(C_{a,w_{l}(x)}\times C_{b,w_{k}(y)})=\\
=\underset{\epsilon\downarrow0}{\lim}\:\frac{\mu((C_{a,w_{l}(x)}\times C_{b,w_{k}(y)})\cap P_{t_{0}}^{-1}(B(P_{t_{0}}z,\epsilon)))}{P_{t_{0}}\mu(B(P_{t_{0}}z,\epsilon))}=\\
=\underset{\epsilon\downarrow0}{\lim}\:\frac{2\epsilon}{P_{t_{0}}\mu(B(P_{t_{0}}z,\epsilon))}\cdot\frac{\mu((C_{a,w_{l}(x)}\times C_{b,w_{k}(y)})\cap P_{t_{0}}^{-1}(B(P_{t_{0}}z,\epsilon)))}{2\epsilon}=\\
=\beta\cdot\underset{\epsilon\downarrow0}{\lim}\:\frac{\mu((C_{a,w_{l}(x)}\times C_{b,w_{k}(y)})\cap P_{t_{0}}^{-1}(B(P_{t_{0}}z,\epsilon)))}{2\epsilon}\:.\label{E8}
\end{multline}
For each $(x',y')\in\mathbb{R}^{2}$ set $g(x',y')=(f_{a,w_{l}(x)}(x'),f_{b,w_{k}(y)}(y'))$,
then
\begin{equation}
C_{a,w_{l}(x)}\times C_{b,w_{k}(y)}=f_{a,w_{l}(x)}(C_{a})\times f_{b,w_{k}(y)}(C_{b})=g(C_{a}\times C_{b}).\label{E9}
\end{equation}
Let $\epsilon>0$, and let $L:\mathbb{R}^{2}\rightarrow\mathbb{R}^{2}$
be a linear map with $L(1,0)=(a^{l},0)$ and $L(0,1)=(0,b^{k})$.
Since $L$ is the linear part of the affine transformation $g$, we
have
\begin{multline}
P_{t_{0}}^{-1}(B(P_{t_{0}}z,\epsilon))=z+V^{t_{0}}+B(0,\epsilon)=\\
=g\circ g^{-1}(z)+L\circ L^{-1}(V^{t_{0}})+L\circ L^{-1}(B(0,\epsilon))=\\
=g(g^{-1}(z)+L^{-1}(V^{t_{0}})+L^{-1}(B(0,\epsilon)))\:.\label{E10}
\end{multline}
From $a^{-l}\ge a^{-t_{0}-k\alpha+1}\geq a\cdot b^{-k}$, we obtain
\begin{equation}
L^{-1}(B(0,\epsilon))\supset B(0,\epsilon\cdot a\cdot b^{-k})\:.\label{E11}
\end{equation}
Also we have
\begin{multline*}
L^{-1}(V^{t_{0}})=L^{-1}((W^{t_{0}})^{\perp})=L^{-1}(((1,\tau\cdot a^{t_{0}})\cdot\mathbb{R})^{\perp})=\\
=L^{-1}((\tau\cdot a^{t_{0}},-1)\cdot\mathbb{R})=(\tau\cdot a^{t_{0}}\cdot a^{-l},-b^{-k})\cdot\mathbb{R}=(\tau\cdot a^{t_{0}}\cdot\frac{b^{k}}{a^{l}},-1)\cdot\mathbb{R},
\end{multline*}
and so since
\[
\frac{b^{k}}{a^{l}}=a^{k\cdot\log_{a}b-l}=a^{k\alpha-[t_{0}+k\alpha]},
\]
it follows that
\begin{equation}
L^{-1}(V^{t_{0}})=(\tau\cdot a^{t_{0}+k\alpha-[t_{0}+k\alpha]},-1)\cdot\mathbb{R}=((1,\tau\cdot a^{R^{k}(t_{0})})\cdot\mathbb{R})^{\perp}=V^{R^{k}(t_{0})}\:.\label{E12}
\end{equation}
Set
\[
Q_{\epsilon}=P_{R^{k}(t_{0})}^{-1}(B(P_{R^{k}(t_{0})}(f_{a,w_{l}(x)}^{-1}(x),f_{b,w_{k}(y)}^{-1}(y)),\epsilon ab^{-k})),
\]
then from (\ref{E10}), (\ref{E11}) and (\ref{E12}) it follows that
\begin{multline*}
P_{t_{0}}^{-1}(B(P_{t_{0}}z,\epsilon))=g(g^{-1}(z)+L^{-1}(V^{t_{0}})+L^{-1}(B(0,\epsilon)))\supset\\
\supset g((f_{a,w_{l}(x)}^{-1}(x),f_{b,w_{k}(y)}^{-1}(y))+V^{R^{k}(t_{0})}+B(0,\epsilon ab^{-k}))=g(Q_{\epsilon}).
\end{multline*}
Now from (\ref{E8}), (\ref{E9}) and Lemma \ref{L9} we get that
\begin{multline*}
\mu_{t_{0},z}(C_{a,w_{l}(x)}\times C_{b,w_{k}(y)})=\beta\cdot\underset{\epsilon\downarrow0}{\lim}\:\frac{\mu(g((C_{a}\times C_{b})\cap Q_{\epsilon}))}{2\epsilon}=\\
=\beta\cdot2^{-l-k}\cdot\frac{a}{b^{k}}\cdot\underset{\epsilon\downarrow0}{\lim}\:\frac{\mu((C_{a}\times C_{b})\cap Q_{\epsilon})}{2\epsilon ab^{-k}}\ge\\
\ge\frac{\beta}{2}\cdot2^{-k-k\alpha}\cdot\frac{a}{b^{k}}\cdot F((f_{a,w_{l}(x)}^{-1}(x),f_{b,w_{k}(y)}^{-1}(y)),R^{k}(t_{0}))=\\
=\frac{\beta a}{2}\cdot2^{-k-k\alpha}\cdot b^{-k}\cdot F(T^{k}(z,t_{0}))\geq\frac{\beta a}{2}\cdot2^{-k-k\alpha}\cdot b^{-k}\cdot M\:.
\end{multline*}
Since 
\[
C_{a,w_{l}(x)}\times C_{b,w_{k}(y)}\subset B(z,\frac{2\cdot b^{k}}{a})\:\mbox{ and }\:2^{-k-k\alpha}\cdot b^{-k}\cdot b^{-k\eta}=1,
\]
it follows that
\begin{multline*}
\frac{\mu_{t_{0},z}(B(z,\frac{2\cdot b^{k}}{a}))}{(4a^{-1}\cdot b^{k})^{\eta}}\geq\frac{\mu_{t_{0},z}(C_{a,w_{l}(x)}\times C_{b,w_{k}(y)})}{(4a^{-1}\cdot b^{k})^{\eta}}\geq\frac{\frac{\beta a}{2}\cdot2^{-k-k\alpha}\cdot b^{-k}\cdot M}{(4a^{-1}\cdot b^{k})^{\eta}}\geq\\
\ge\frac{\beta a^{2}}{8}\cdot M\cdot2^{-k-k\alpha}\cdot b^{-k}\cdot b^{-k\eta}=\frac{\beta a^{2}}{8}\cdot M\:.
\end{multline*}
This shows that $\Theta^{*\eta}(\mu_{t_{0},z},z)\geq\frac{\beta a^{2}}{8}\cdot M$,
which proves (\ref{E7}) since $\beta>0$ and $M$ can be chosen arbitrarily
large.

Let $z\in A$ and $u\in A\cap V_{z}^{t_{0}}$, then from (\ref{E7})
\[
\Theta^{*\eta}(\mu_{t_{0},z},u)=\Theta^{*\eta}(\mu_{t_{0},u},u)=\infty,
\]
and so from Theorem 6.9 in \cite{key-1} we get that $\mathcal{H}^{\eta}(A\cap V_{z}^{t_{0}})=0$.
Also it holds that $\mu(K\setminus A)=0$, hence from Theorem 7.7
in \cite{key-1} and from Lemma \ref{L10} we get that
\[
\int_{W^{t_{0}}}\mathcal{H}^{\eta}((K\setminus A)\cap V_{u}^{t_{0}})\: d\mathcal{H}^{1}(u)\leq const\cdot\mathcal{H}^{\eta+1}(K\setminus A)=const\cdot\mu(K\setminus A)=0.
\]
This shows that $\mathcal{H}^{\eta}((K\setminus A)\cap V_{u}^{t_{0}})=0$
for $\mathcal{H}^{1}$-a.e. $u\in W^{t_{0}}$, and so $\mathcal{H}^{\eta}((K\setminus A)\cap V_{z}^{t_{0}})=0$
for $\mu$-a.e. $z\in K$ since $P_{t_{0}}\mu\ll\mathcal{H}^{1}$.
It follows that for $\mu$-a.e. $z\in A$, and so for $\mu$-a.e.
$z\in K$,
\[
\mathcal{H}^{\eta}(K\cap V_{z}^{t_{0}})=\mathcal{H}^{\eta}(A\cap V_{z}^{t_{0}})+\mathcal{H}^{\eta}((K\setminus A)\cap V_{z}^{t_{0}})=0\:.
\]
From this, from Lemma \ref{L4}, and from Fubini's theorem it follows
that $\mathcal{H}^{\eta}(K\cap V_{z}^{t})=0$ for $\nu$-a.e. $(z,t)\in X$,
which completes the proof of Theorem \ref{T13}. $\square$$\newline$

\emph{Proof of Theorem \ref{T6}:} Let $G$ be the set of all $1$-dimensional
linear subspaces of $\mathbb{R}^{2}$, and set
\[
E=\{(z,V)\in K\times G\::\:\mathcal{H}^{d_{a}+d_{b}-1}(K\cap V_{z})=0\}.
\]
For each $-\infty\le t_{1}<t_{2}\le\infty$ set
\[
G_{t_{1},t_{2}}=\{V\in G\::\: V=(t,-1)\cdot\mathbb{R}\mbox{ with }t\in(t_{1},t_{2})\}.
\]
Given $k\in\mathbb{Z}$ we can apply theorem \ref{T13} with $\tau=a^{k}$,
in order to get that $(z,V)\in E$ for $\mu\times\xi_{G}$-a.e. $(z,V)\in K\times G_{a^{k+1},a^{k}}$.
By doing this for each $k\in\mathbb{Z}$ we get that $(z,V)\in E$
for $\mu\times\xi_{G}$-a.e. $(z,V)\in K\times G_{0,\infty}$. Now
Theorem \emph{\ref{T6}} follows from the symmetry of $K$ with respect
to the map that takes $(x,y)\in K$ to $(1-x,y)$. $\square$

\end{document}